\documentclass[11pt]{article}

\usepackage{amssymb,amsmath, amsfonts, amsthm}

\usepackage{graphicx}
\usepackage{pgfplots}
\usepackage{tikz}
\usepackage{caption}
\usepackage{booktabs}
\usepackage{array}
\usepackage{verbatim}
\usepackage{subfig}
\usepackage{geometry}
\usepackage{fancyhdr}

\newcommand{\vertex}{\node[vertex]}
\tikzstyle{vertex}=[circle, draw, inner sep=0pt, minimum size=6pt]
\newtheorem{theorem}{Theorem}
\newtheorem{lemma}{Lemma}

\newtheorem{obs}{Observation}

\newcommand{\smallqed}{{\tiny ($\Box$)}}
\newcommand{\dist}{\hskip2pt$\rm{dist}$}

\begin{document}

\title{Characterizing all $K_4$-free well-edge-dominated graphs of girth $3$}

\author{Sarah E. Anderson$^a$ \and Kirsti Kuenzel$^b$
}

\maketitle

\begin{center}
$^a$ Department of Mathematics, University of St. Thomas, St. Paul, Minnesota, USA\\
$^b$ Department of Mathematics, Trinity College, Hartford, CT, USA\\

\end{center}
\medskip

\maketitle
\begin{abstract}
Given a graph $G$, a set $F$ of edges is an edge dominating set if all edges in $G$ are either in $F$ or adjacent to an edge in $F$. $G$ is said to be well-edge-dominated if every minimal edge dominating set is also minimum. In 2022, it was proven that there are precisely three nonbipartite, well-edge-dominated graphs with girth at least four. Then in 2025, a characterization of all well-edge-dominated graphs containing exactly one triangle was found. In this paper, we characterize all well-edge-dominated graphs that contain a triangle and yet are $K_4$-free.
\end{abstract}

%%%%%%%%%%%%%%%%%%%%%%%%%%%%%%%%%%%
%%%%%%%%%%%%%%%%%%%%%%%%%%%%%%%%%%%
\maketitle
\section{Introduction}
Given a graph $G$, a set of edges $F$ in $G$ is called a \emph{matching} if no pair of edges in $F$ share a common vertex, and it is called an edge dominating set if each edge in $G$ is adjacent to an edge in $F$. We typically want to maximize the cardinality of a matching, giving us the matching number $\alpha'(G)$, and  minimize the cardinality of an edge dominating set, giving us the edge domination number $\gamma'(G)$. $G$ is said to be \emph{equimatchable} if all maximal matchings in $G$ have the same cardinality. This property was first defined, independently, by Lewin \cite{L-1974} and Meng \cite{M-1974}. Later  Lesk, Plummer, and Pulleyblank \cite{LPP-1984} presented a polynomial time algorithm for recognizing when a graph is equimatchable, and since then, a search for a structural characterization has been a running theme. $G$ is said to be \emph{well-edge-dominated} if all minimal edge dominating sets in $G$ have the same cardinality. Seeing as how any maximal matching in $G$ is also a minimal edge dominating set, it is no surprise that if $G$ is well-edge-dominated, then it is also equimatchable. 

The class of equimatchable graphs are well-studied (see \cite{AGHI-2018, BGO-2020, EK-2016, KPS-2003}). Notably, Frendrup, Hartnell, and Vestergaard \cite{FHV-2010} proved that a connected equimatchable graph with no cycles of length less than $5$ is either $C_5$, $C_7$, or belongs to the family of bipartite graphs $\mathcal{C}$ that contains $K_2$ as well as those with a bipartition $A\cup B$ where $A$ is the set of support vertices. Then in \cite{BGO-2020}, B\"{u}y\"{u}k\c{c}olak et al. characterized all nonbipartite equimatchable graphs containing $4$-cycles.  It is worth mentioning that although the equimatchable graphs were characterized in \cite{LPP-1984}, no structural results were provided so it was unclear what the class of bipartite graphs containing $4$-cycles look like. Since then, several papers have attempted to give a structural characterization of the bipartite equimatchable graphs containing $4$-cycles (see \cite{BGO-2023, DE-2019}). Frendrup, Hartnell, and Vestergaard \cite{FHV-2010} also proved that an equimatchable graph with girth at least $5$ is well-edge-dominated. In \cite{AKR-2022} Anderson et al. showed that there are only three nonbipartite well-edge-dominated graphs with girth at least $4$; namely, $C_5$, $C_7$, and $C_7^*$ which is obtained from $C_7$ by adding a chord to $C_7$ between two vertices at distance $3$. Then in \cite{BCKKPRV-2025}, Berg et al.  characterized all well-edge-dominated graphs containing exactly one triangle. Based on their construction, we characterize all $K_4$-free well-edge-dominated graphs of girth $3$. 

The remainder of the paper is organized as follows. In Section~\ref{sec:Previous} we provide definitions and known results that we rely on in later sections. In Section~\ref{sec:K_4-free} we provide an infinite class of nonbipartite $K_4$-free well-edge-dominated graphs which we refer to as $\mathcal{G}$, and we prove that any nonbipartite $K_4$-free well-edge-dominated graph is either in $\mathcal{G}$ or one of six exceptions.

  \subsection{Previous Results}\label{sec:Previous}
Let $G= (V(G), E(G))$ be any graph where we write $n(G)$ to denote $|V(G)|$. Given any $v\in V(G)$, the open neighborhood of $v$ is the set $N_G(v) = \{u \in V(G): uv \in E(G)\}$ (or simply $N(v)$ when the context is clear). The closed neighborhood is defined to be $N_G[v] = N_G(v) \cup \{v\}$ (or simply $N[v]$ when the context is clear). For any set $S \subseteq V(G)$, $N_G(S) = \cup_{v \in S} N_G(v)$. Similarly, given any edge $f \in E(G)$, the open edge neighborhood of $f$, denoted $N_e(f)$, is the set of all edges in $G$ adjacent to $f$, and the closed edge neighborhood of $f$ is $N_e[f]= N_e(f) \cup\{f\}$. For any set $F\subseteq E(G)$, the closed edge neighborhood of $F$ is $N_e[F] = \cup_{f \in F} N_e[f]$. For any $v \in V(G)$, the degree of $v$ is $\deg_G(v) = |N_G(v)|$ (or simply $\deg(v)$). Any vertex with degree one is called a leaf, and its lone neighbor in $G$ is referred to as a support vertex in $G$. Given any two vertices $u,v\in V(G)$, the graph obtained by \emph{identifying $u$ and $v$} is found by removing $u$ and $v$ from $G$ and replacing with a vertex $w$ where $w$ is adjacent to each vertex in $(N_G(u) \cup N_G(v))- \{u,v\}$. We will refer to the operation of \emph{identifying each vertex in $S \subset V(G)$ into a single vertex} to mean removing each vertex in $S$ and replacing with a vertex $w$ where $w$ is adjacent to each vertex in $N_G(S) - S$. Given a graph $F$, $G$ is said to $F$-free if $G$ does not contain $F$ as an induced subgraph. 

An edge dominating set is a set $D \subseteq E(G)$ such that every edge in $G$ is either in $D$ or adjacent to an edge in $D$. The edge domination number of $G$, denoted $\gamma'(G)$, is the minimum cardinality among all edge dominating sets of $G$. Given $e \in D$, the private edge neighbors of $e$ with respect to $D$ is the set of edges $f$ where $N_e[f] \cap D =\{e\}$.  A matching in $G$ is a set of independent edges in $G$, and the matching number, denoted $\alpha'(G)$, is the maximum cardinality of a matching in $G$. For any graph $G$, $\gamma'(G) \le \alpha'(G)$. Recall that $G$ is bipartite if we can partition $V(G)$ into two independent sets. If $G$ is bipartite, we write $G = (A\cup B, E)$ to mean that $A$ and $B$ are the partite sets in $G$. 

$G$ is said to be equimatchable if every maximal matching in $G$ has the same cardinality, namely, $\alpha'(G)$.  Seeing as how our constructions for well-edge-dominated graphs rely on equimatchable bipartite graphs, one important  result regarding equimatchable bipartite  graphs is the following. 

  \begin{lemma}\cite{BGO-2023} Let $G = (A\cup B, E)$ with $|A| < |B|$ be a connected equimatchable bipartite graph. Then each vertex $u \in A$ satisfies at least one of the following.
  \begin{enumerate}
  \item[(i)] $u$ is a support vertex  in $G$,
  \item[(ii)] $u$ is included in a subgraph $K_{2, 2}$ in $G$.
  \end{enumerate}
  \end{lemma}
  
  The next question would be what can be said about bipartite equimatchable graphs $G = (A \cup B, E)$ where $|A| = |B|$?   A graph $G$ is defined to \emph{randomly matchable} if it is an equimatchable graph admitting a perfect matching. Sumner characterized all randomly matchable graphs. 
  
  \begin{theorem}\cite{S-1979}\label{thm:rmatch} A connected graph is randomly matchable if and only if it is isomorphic to $K_{2n}$ or $K_{n,n}$ for $n \ge 1$. 
  \end{theorem}
  
 In this paper, we focus on graphs that are well-edge-dominated, i.e. graphs where every minimal edge dominating set is also minimum. There is a well-known chain of inequalities, namely 
 \[\gamma'(G) \le \alpha'(G) \le \Gamma'(G)\]
where  $\Gamma'(G)$ is the maximum cardinality among all minimal dominating sets. Thus, if $G$ is well-edge-dominated, then $\gamma'(G) = \Gamma'(G)$, implying that $G$ is equimatchable. This together with the fact that all  equimatchable graphs of girth $5$ or more are characterized meant that the only well-edge-dominated graphs left to classify were those of girth $3$ or $4$.  In 2022, Anderson et al. \cite{AKR-2022} proved that there are only three nonbipartite well-edge-dominated graphs of girth $4$. We let $C_7^*$ denote the graph obtained from $C_7$ by adding a chord between one pair of vertices at distance three. 
  \begin{theorem}\label{thm:girth5red}\cite{AKR-2022} If $G$ is well-edge-dominated with $g(G) \ge 4$, then either $G$ is bipartite or $G \in \{C_5, C_7, C_7^*\}$.
  \end{theorem}
  
\begin{figure}[h!]
\begin{center}
\begin{tikzpicture}[scale=.65]
    \vertex (1) at (0,0)  [fill, scale=.75, label=below:$a$]{};
    \vertex (2) at (0,2) [fill, scale=.75, label=left:$c$]{};
    \vertex (3) at (3,0)  [fill, scale=.75, label=below:$b$]{};
    \vertex (4) at (3,2)   [fill, scale=.75, label=right:$d$]{};
    \vertex (5) at (0,4) [fill, scale=.75, label=left:$x$]{};
    \vertex (6) at (3,4)  [fill, scale=.75, label=right:$y$]{};
    \vertex (7) at (1.5, 5)  [fill, scale=.75, label=above:$z$]{};

    \vertex (8) at (6.5,0) [fill, scale=.75]{};
    \vertex (9) at (9.5,0)  [fill, scale=.75]{};
    \vertex (10) at (6.5,2)  [fill, scale=.75]{};
    \vertex (11) at (9.5,2)  [fill, scale=.75]{};
    \vertex (12) at (8,3)  [fill, scale=.75]{};
    
     \vertex (1a) at (14,3)  [scale=.75, fill=black, label=above:$a$]{};
    \vertex (2a) at (12.5,1.5)  [scale=.75, fill=black, label=left:$d$]{};
    \vertex (3a) at (14,0)  [scale=.75, fill=black, label=below:$c$]{};
    \vertex (4a) at (15.5,1.5)  [scale=.75, fill=black, label=above:$b$]{};
    \vertex (5a) at (17.5,1.5)  [scale=.75, fill=black, label=above:$y$]{};
    \vertex (6a) at (19,3)  [scale=.75, fill=black, label=above:$x$]{};
    \vertex (7a) at (19,0)  [scale=.75, fill=black, label=below:$z$]{};
    \node(A) at (1.4, -1)[]{(a) $\mathcal{DH}$};
    \node(B) at (8, -1)[]{(b) $\mathcal{H}$};
    \node(C) at (16.5, -1)[]{(c) $Cr$};

    \path 
    (1) edge (2)
    (1) edge (3)
    (4) edge (2)
    (4) edge (3)
    (5) edge (2)
    (6) edge (4)
    (7) edge (6)
    (7) edge (5)
    (5) edge (6)

    (8) edge (9)
    (8) edge (10)
    (11) edge (10)
    (11) edge (9)
    (11) edge (12)
    (12) edge (10)
    
       (1a) edge (2a)
    (2a) edge (3a)
    (3a) edge (4a)
    (4a) edge (1a)
    (4a) edge (5a)
    (5a) edge (6a)
    (5a) edge (7a)
    (6a) edge (7a)
    (6a) edge (1a)
    (7a) edge (3a)
    ;

\end{tikzpicture}
\caption{The dream house (a), the house graph (b) and the crystal graph (c)}
\label{fig:houses}
\end{center}
\end{figure}
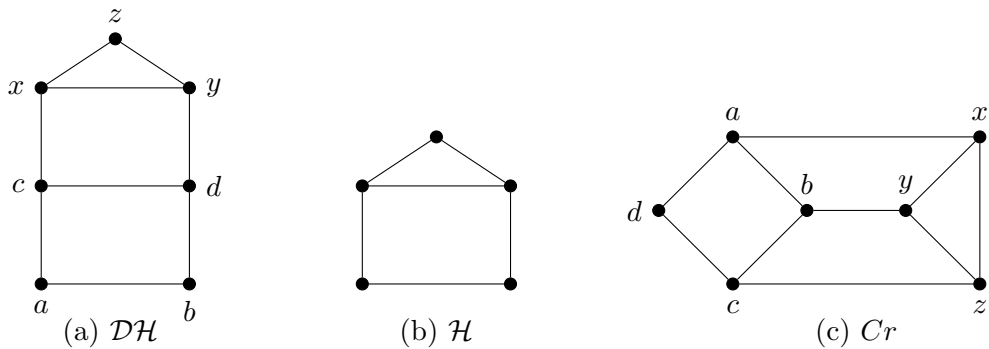

Then in \cite{BCKKPRV-2025}, a characterization of well-edge-dominated graphs containing exactly one triangle was provided. It relied on two infinite families of graphs, each of which were constructed from well-edge-dominated bipartite graphs (of which we still do not have a nice structural description, but nevertheless). Given a well-edge-dominated bipartite graph $G= (A \cup B, E)$  such that $|A| < |B|$, we say that $v \in B$ is \emph{detachable} if $G-v$ is well-edge-dominated. Note that in this case $\gamma'(G-v) = |A|$ as $G-v$ contains a matching of size $|A|$. Additionally, we say that $v \in B$ is \emph{strongly detachable} if $G-v$ is well-edge-dominated and each vertex in $N_G(v)$ is a support vertex  in $G-v$. 
The first class of well-edge-dominated graphs containing exactly one triangle is denoted $\mathcal{F}$ and defined as follows. We say that $G \in \mathcal{F}$ if $G$ is obtained from the disjoint union of the house graph $\mathcal{H}$, depicted in Figure~\ref{fig:houses}(b), and a well-edge-dominated bipartite graph $G'=( A \cup B, E)$ where  $|A|< |B|$,  by identifying the vertex of degree two in $\mathcal{H}$ on the triangle with $w \in V(G')$ where $w$ is strongly detachable.  Examples of such graphs are $W_7$ and $V_7$, found in Figures~\ref{ND7} and \ref{ND8}, respectively. The second class of well-edge-dominated graphs containing exactly one triangle is denoted $\mathcal{T}$ and defined as follows. We say that $G \in \mathcal{T}$ if $G$ is obtained from the disjoint union of $K_3$ and a well-edge-dominated bipartite graph $G'=(A \cup B, E)$ where  $|A| < |B|$,  by identifying $x \in V(K_3)$ with $w \in V(G')$ where $w$ is detachable. Examples of such graphs are $W_5, W_6, W_9, V_4, V_5$, and $V_6$ in Figures~\ref{ND7} and \ref{ND8}. Letting $Cr$ denote the ``crystal" graph depicted in Figure~\ref{fig:houses}(c) and $\mathcal{DH}$ denote the ``dream house" graph depicted in Figure~\ref{fig:houses}(a), the complete class of well-edge-dominated graphs containing exactly one triangle is the following.

 \begin{theorem}\cite{BCKKPRV-2025}\label{thm:onetriangle} $G$ is a well-edge-dominated graph with exactly one triangle if and only if $G \in \mathcal{T}\cup \mathcal{F} \cup \{K_3, Cr, \mathcal{H}, \mathcal{DH}\}$.
 \end{theorem}

  %Therefore, if $G$ is an equimatchable bipartite graph $G = (U\cup V, E)$ where $|U| = |V|$, then it contains a maximal matching that saturates $U$ and therefore $V$, so it admits a perfect matching. This implies that $G = K_{n,n}$ for $n \ge 1$. 

In this paper, we characterize all $K_4$-free well-edge-dominated graphs containing a $3$-cycle. In that vein, we will need the following preliminary results used in the next section. 

\begin{lemma}\label{lem:noleaves}\cite{BCKKPRV-2025} If $G$ is a graph of order at least $4$ obtained from the triangle $xyz$ by appending leaves to $x, y$, or $z$, then $G$ is not well-edge-dominated. 
\end{lemma}

Additionally, the ``edge version'' of a fact used by Finbow,
Hartnell and Nowakowski in~\cite{fhn-1988} is that provided $G$ is well-edge-dominated, then $M \cup D_1$ and $M \cup D_2$ are both minimal edge dominating sets of $G$
for any matching $M$ and any pair $D_1$ and $D_2$ of minimal edge dominating sets of the graph $G-N_e[M]$. This was first observed in \cite{AKR-2022}.

\begin{lemma}\label{lem:reduce}\cite{AKR-2022} Let $M$ be a matching in a graph $G$. If $G$ is well-edge-dominated, then $G - N_e[M]$ is well-edge-dominated. If $G$ is equimatchable, then $G - N_e[M]$ is equimatchable. 
\end{lemma}

We also need the following result regarding well-edge-dominated bipartite graphs. 

\begin{lemma}\label{lem:bipartitesupports} Let $G = (A\cup B, E)$ be a well-edge-dominated bipartite graph where $|A| < |B|$. If $y \in A$ is not a support vertex in $G$, then there exists a minimal edge dominating set $F$ of $G$ that does not contain an edge incident to $y$.
\end{lemma}

\begin{proof} To see this, we construct such an edge dominating set $F$. First, for each $b \in N_G(y)$, choose an edge incident to $b$ and a vertex in $A- \{y\}$ and call the resulting set $F_0$. Let $A' = \{a_1, \dots, a_k\}$ be the vertices in $A- \{y\}$ that are not incident to an edge in $F_0$. We construct sets $F_i$ for each $i \in [k]$ in the following way.  If there exists an edge $a_1b'$ such that $b'$ is not incident to an edge in $F_0$, set $F_1 = F_0 \cup \{a_1b'\}$. Otherwise, set $F_1 = F_0$. Continue this process for each $2 \le i \le k$ until we arrive at the set $F_k$. We claim that $F_k$ is a minimal edge dominating set of $G$. Note that if $x \in A - \{y\}$ is not incident to an edge in $F_k$, then every neighbor of $x$ in $B$ is incident to an edge in $F_k$ by construction. So $F_k$ is indeed an edge dominating set of $G$. Moreover, each edge $e \in F_k - F_0$ is its own private edge neighbor with respect to $F_k$. Finally, each edge in $F_0$ is of the form $bx$ where $x \in A- \{y\}$ and $b \in N_G(y)$. Thus, the private edge neighbor of $bx$ is $by$. 
\end{proof}

Finally, we introduce two other graphs that will show up throughout Section~\ref{sec:K_4-free}. The graph $F_5$ is the fan on $5$ vertices. A \emph{diamond} is the graph $K_4 - e$ for any edge $e$. Given a diamond with vertices $\{w, x, y, z\}$ where $w$ and $z$ have degree $2$, we refer to $w$ and $z$ as the \emph{exterior vertices} of the diamond, and we refer to $x$ and $y$ as the \emph{interior vertices} of the diamond.

%%%%%%%%%%%%%%%%%%%%%%%%%%%%%%%%%%%%%%%%%%%%%%%%%%%%%%%%%%%%%%%%%%%%%%%%%%%%%%%%%%
\section{$K_4$-free well-edge-dominated graphs}\label{sec:K_4-free}
\begin{figure}[h!]
\begin{center}
\begin{tikzpicture}[scale=.85]

    \vertex (01) at (0.5,7)  [scale=.75, fill=black, label=below:$$]{};
    \vertex (02) at (0.5,8)  [scale=.75, fill=black, label=left:$$]{};
    \vertex (03) at (1.5, 7.5)  [scale=.75, fill=black]{};
    \vertex (04) at (1.5, 8.5)  [scale=.75, fill=black, label=above:$$]{};
    \vertex (05) at (2.5, 8)  [scale=.75, fill=black, label=right:$$]{};
    \vertex (06) at (2.5, 7)  [scale=.75, fill=black, label=below:$$]{};
    \vertex (07) at (1.5, 6.5)  [scale=.75, fill=black, label=below:$$]{};
  
     \vertex (0B1) at (5,7)  [scale=.75, fill=black, label=below:$$]{};
    \vertex (0B2) at (5,8)  [scale=.75, fill=black, label=left:$$]{};
    \vertex (0B3) at (6, 7.5)  [scale=.75, fill=black]{};
    \vertex (0B4) at (6, 8.5)  [scale=.75, fill=black, label=above:$$]{};
    \vertex (0B5) at (7, 8)  [scale=.75, fill=black, label=right:$$]{};
    \vertex (0B6) at (7, 7)  [scale=.75, fill=black, label=below:$$]{};
    \vertex (0B7) at (6, 6.5)  [scale=.75, fill=black, label=below:$$]{};
    
       \vertex (0C1) at (9.5,7)  [scale=.75, fill=black, label=below:$$]{};
    \vertex (0C2) at (9.5,8)  [scale=.75, fill=black, label=left:$$]{};
    \vertex (0C3) at (10.5, 7.5)  [scale=.75, fill=black]{};
    \vertex (0C4) at (10.5, 8.5)  [scale=.75, fill=black, label=above:$$]{};
    \vertex (0C5) at (11.5, 8)  [scale=.75, fill=black, label=right:$$]{};
    \vertex (0C6) at (11.5, 7)  [scale=.75, fill=black, label=below:$$]{};
    \vertex (0C7) at (10.5, 6.5)  [scale=.75, fill=black, label=below:$$]{};

       \vertex (0D1) at (14,7)  [scale=.75, fill=black, label=below:$$]{};
    \vertex (0D2) at (14,8)  [scale=.75, fill=black, label=left:$$]{};
    \vertex (0D3) at (15, 7.5)  [scale=.75, fill=black]{};
    \vertex (0D4) at (15, 8.5)  [scale=.75, fill=black, label=above:$$]{};
    \vertex (0D5) at (16, 8)  [scale=.75, fill=black, label=right:$$]{};
    \vertex (0D6) at (16, 7)  [scale=.75, fill=black, label=below:$$]{};
    \vertex (0D7) at (15, 6.5)  [scale=.75, fill=black, label=below:$$]{};

    \node(A) at (1.5, 5) []{$W_1$};
    \node(B) at (6, 5)[]{$W_2$};
    \node(C) at (10.5, 5)[]{$W_3$};
    \node(D) at (15, 5)[]{$W_4$};

    \vertex (1) at (.5,2.5)  [scale=.75, fill=black, label=below:]{};
    \vertex (2) at (1.5, 2.5)  [scale=.75, fill=black, label=left:]{};
    \vertex (3) at (1, 1.5)  [scale=.75, fill=black]{};
    \vertex (4) at (1, .5)  [scale=.75, fill=black, label=above:]{};
    \vertex (5) at (.5, -.5)  [scale=.75, fill=black, label=right:]{};
    \vertex (6) at (1, -.5)  [scale=.75, fill=black, label=below:]{};
    \vertex (7) at (1.5, -.5)  [scale=.75, fill=black, label=below:]{};
  
     \vertex (B1) at (4,2.5)  [scale=.75, fill=black, label=below:]{};
    \vertex (B2) at (5,2.5)  [scale=.75, fill=black, label=left:]{};
    \vertex (B3) at (4.5, 1.5)  [scale=.75, fill=black]{};
    \vertex (B4) at (4, .5)  [scale=.75, fill=black, label=above:]{};
    \vertex (B5) at (5, .5)  [scale=.75, fill=black, label=right:]{};
    \vertex (B6) at (4, -.5)  [scale=.75, fill=black, label=below:]{};
    \vertex (B7) at (5, -.5)  [scale=.75, fill=black, label=below:]{};

        \vertex (C1) at (8,2.5)  [scale=.75, fill=black, label=below:]{};
    \vertex (C2) at (9,2.5)  [scale=.75, fill=black, label=left:]{};
    \vertex (C3) at (8, 1.5)  [scale=.75, fill=black]{};
    \vertex (C4) at (9, 1.5)  [scale=.75, fill=black, label=above:]{};
    \vertex (C5) at (8.5, .5)  [scale=.75, fill=black, label=right:]{};
    \vertex (C6) at (8.5, 0)  [scale=.75, fill=black, label=below:]{};
    \vertex (C7) at (8.5, -.5)  [scale=.75, fill=black, label=below:]{};

        \vertex (D1) at (11.5,2.5)  [scale=.75, fill=black, label=below:]{};
    \vertex (D2) at (12.5,2.5)  [scale=.75, fill=black, label=left:]{};
    \vertex (D3) at (11.5, 1.5)  [scale=.75, fill=black]{};
    \vertex (D4) at (12.5, 1.5)  [scale=.75, fill=black, label=above:]{};
    \vertex (D5) at (12, .5)  [scale=.75, fill=black, label=right:]{};
    \vertex (D6) at (11.5, -.5)  [scale=.75, fill=black, label=below:]{};
    \vertex (D7) at (12.5, -.5)  [scale=.75, fill=black, label=below:]{};
    
     \vertex (E1) at (15,2.5)  [scale=.75, fill=black, label=below:]{};
    \vertex (E2) at (16,2.5)  [scale=.75, fill=black, label=left:]{};
    \vertex (E3) at (15, 1.5)  [scale=.75, fill=black]{};
    \vertex (E4) at (16, 1.5)  [scale=.75, fill=black, label=above:]{};
    \vertex (E5) at (15.5, .5)  [scale=.75, fill=black, label=right:]{};
    \vertex (E6) at (15, -.5)  [scale=.75, fill=black, label=below:]{};
    \vertex (E7) at (16, -.5)  [scale=.75, fill=black, label=below:]{};
    
         \vertex (F1) at (1.5,-3.5)  [scale=.75, fill=black, label=below:]{};
    \vertex (F2) at (2.5,-3.5)  [scale=.75, fill=black, label=left:]{};
    \vertex (F3) at (1.5, -4.5)  [scale=.75, fill=black]{};
    \vertex (F4) at (2.5, -4.5)  [scale=.75, fill=black, label=above:]{};
    \vertex (F5) at (1.5, -5.5)  [scale=.75, fill=black, label=right:]{};
    \vertex (F6) at (2.5, -5.5)  [scale=.75, fill=black, label=below:]{};
    \vertex (F7) at (2, -6.5)  [scale=.75, fill=black, label=below:]{};

         \vertex (G1) at (6,-3.5)  [scale=.75, fill=black, label=below:]{};
    \vertex (G2) at (7,-3.5)  [scale=.75, fill=black, label=left:]{};
    \vertex (G3) at (6.5, -4.5)  [scale=.75, fill=black]{};
    \vertex (G4) at (6, -5.5)  [scale=.75, fill=black, label=above:]{};
    \vertex (G5) at (7, -5.5)  [scale=.75, fill=black, label=right:]{};
    \vertex (G6) at (6.5, -6)  [scale=.75, fill=black, label=below:]{};
    \vertex (G7) at (6.5, -6.5)  [scale=.75, fill=black, label=below:]{};

         \vertex (H1) at (10,-3.5)  [scale=.75, fill=black, label=below:]{};
    \vertex (H2) at (11,-3.5)  [scale=.75, fill=black, label=left:]{};
    \vertex (H3) at (10.5, -4.5)  [scale=.75, fill=black]{};
    \vertex (H4) at (10.5, -5.5)  [scale=.75, fill=black, label=above:]{};
    \vertex (H5) at (10, -6)  [scale=.75, fill=black, label=right:]{};
    \vertex (H6) at (11, -6)  [scale=.75, fill=black, label=below:]{};
    \vertex (H7) at (10.5, -6.5)  [scale=.75, fill=black, label=below:]{};
    
      \vertex (I1) at (14.25,-4)  [scale=.75, fill=black, label=below:]{};
    \vertex (I2) at (15.25,-4)  [scale=.75, fill=black, label=left:]{};
    \vertex (I3) at (14.75, -5)  [scale=.75, fill=black]{};
    \vertex (I4) at (13.75, -5.25)  [scale=.75, fill=black, label=above:]{};
    \vertex (I5) at (14.25, -6)  [scale=.75, fill=black, label=right:]{};
    \vertex (I6) at (15.25, -6)  [scale=.75, fill=black, label=below:]{};
    \vertex (I7) at (15.75, -5.25)  [scale=.75, fill=black, label=below:]{};

    \node(A) at (1, -1.5) []{$W_5$};
    \node(B) at (4.5, -1.5)[]{$W_6$};
    \node(C) at (8.5, -1.5)[]{$W_7$};
     \node(D) at (12, -1.5)[]{$W_8$};
      \node(E) at (15.5, -1.5)[]{$W_9$};
       \node(F) at (2, -7.5)[]{$W_{10}$};
         \node(G) at (6.5, -7.5)[]{$W_{11}$};
           \node(H) at (10.5, -7.5)[]{$W_{12}$};
             \node(I) at (14.5, -7.5)[]{$W_{13}$};

    \path 
    	(01) edge (02)
	(02) edge (03)
	(03) edge (04)
	(04) edge (05)
	(05) edge (06)
	(06) edge (07)
	(01) edge (07)
	(02) edge (04)
	(03) edge (05)
	
	(0B1) edge (0B2)
	(0B2) edge (0B3)
	(0B3) edge (0B4)
	(0B4) edge (0B5)
	(0B5) edge (0B6)
	(0B6) edge (0B7)
	(0B1) edge (0B7)
	(0B2) edge (0B4)
	(0B3) edge (0B5)
	(0B3) edge (0B7)
	
  	(0C1) edge (0C2)
	(0C2) edge (0C3)
	(0C3) edge (0C4)
	(0C4) edge (0C5)
	(0C5) edge (0C6)
	(0C6) edge (0C7)
	(0C1) edge (0C7)
	(0C2) edge (0C4)
	(0C3) edge (0C5)  
	(0C4) edge (0C6)
	
	(0D1) edge (0D2)
	(0D2) edge [bend right] (0D6)
	(0D2) edge  (0D7)
	(0D2) edge (0D3)
	(0D3) edge (0D4)
	(0D4) edge (0D5)
	(0D2) edge (0D4)
	(0D3) edge (0D5)  
	
	(1) edge (2)
	(1) edge (3)
	(2) edge (3)
	(4) edge (3)
	(4) edge (5)
	(4) edge (6)
	(4) edge (7)
	
	(B1) edge (B2)
	(B1) edge (B3)
	(B2) edge (B3)
	(B3) edge (B4)
	(B3) edge (B5)
	(B4) edge (B6)
	(B5) edge (B7)

	(C1) edge (C2)
  	(C1) edge (C3)
	(C2) edge (C4)
	(C3) edge (C4)
	(C4) edge (C5)
	(C3) edge (C5)
	(C5) edge (C6)
	(C6) edge (C7)
	
		(D1) edge (D2)
  	(D1) edge (D3)
	(D2) edge (D4)
	(D3) edge (D4)
	(D4) edge (D5)
	(D3) edge (D5)
	(D5) edge (D6)
	(D5) edge (D7)
	(D6) edge (D7)
	
		(E1) edge (E2)
  	(E1) edge (E3)
	(E2) edge (E4)
	(E3) edge (E4)
	(E4) edge (E5)
	(E5) edge (E6)
	(E5) edge (E7)
	(E6) edge (E7)
	
		(F1) edge (F2)
  	(F1) edge (F3)
	(F2) edge (F4)
	(F3) edge (F4)
	(F4) edge (F6)
		(F3) edge (F5)
	(F5) edge (F6)
	(F5) edge (F7)
	(F6) edge (F7)
	
	(G1) edge (G2)
  	(G1) edge (G3)
	(G2) edge (G3)
	(G3) edge (G4)
	(G3) edge (G5)
		(G4) edge (G5)
	(G3) edge (G6)
	(G6) edge (G7)

	(H1) edge (H2)
  	(H1) edge (H3)
	(H2) edge (H3)
	(H3) edge (H4)
	(H4) edge (H5)
		(H5) edge (H6)
		(H1) edge (H5)
  	(H2) edge (H6)
	(H4) edge (H3)
	(H7) edge (H4)
	(H6) edge (H7)
	
	(I1) edge (I2)
  	(I1) edge (I3)
	(I2) edge (I3)
	(I3) edge (I4)
	(I4) edge (I5)
		(I5) edge (I3)
  	(I3) edge (I6)
	(I3) edge (I7)
	(I6) edge (I7)
	
    ;

\end{tikzpicture}
\caption{All connected, $K_4$-free, well-edge dominated graphs containing a triangle of order $7$ where only $W_1, W_2, W_3$, and $W_4$ contain a diamond}
\label{ND7}
\end{center}
\end{figure}
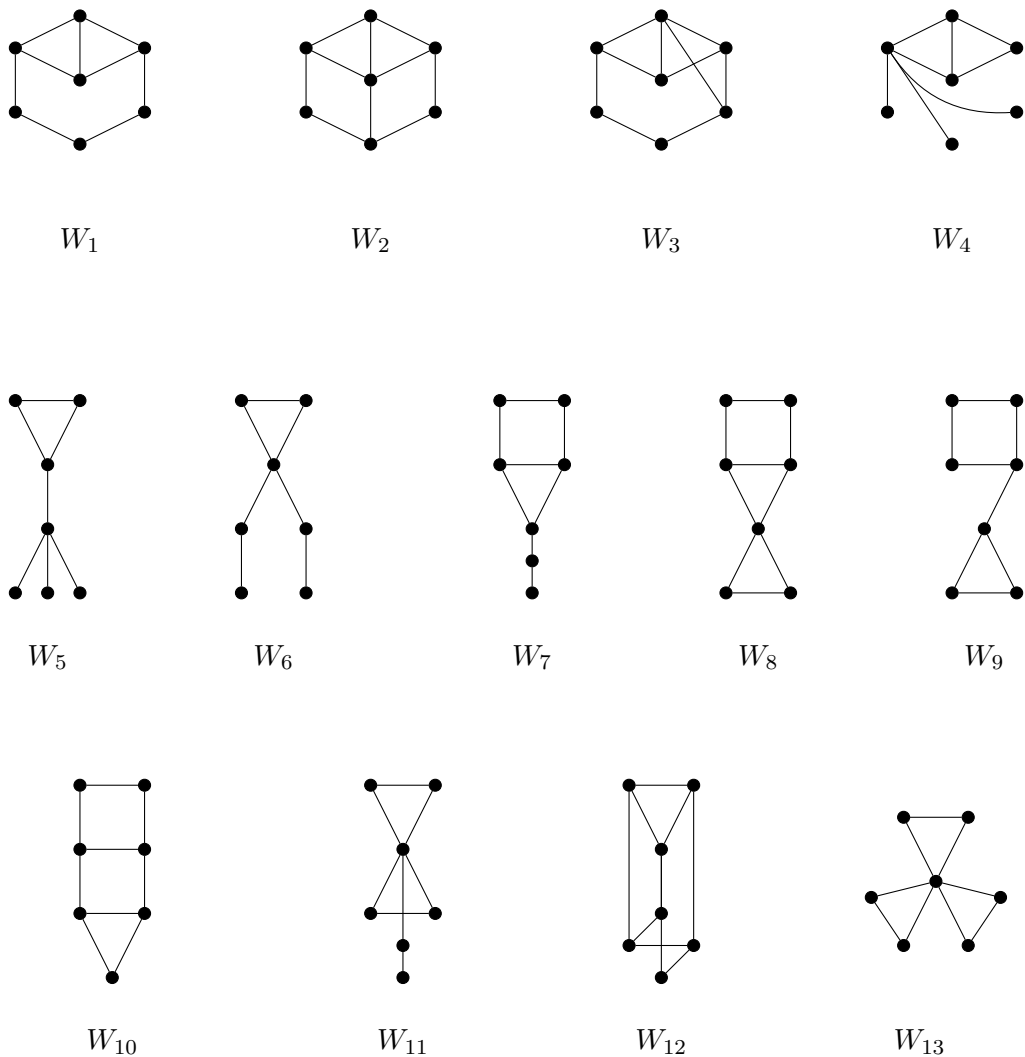

\begin{figure}[h!]
\begin{center}
\begin{tikzpicture}[scale=.85]

    \vertex (01) at (0,7)  [scale=.75, fill=black, label=below:]{};
    \vertex (02) at (.5,8.5)  [scale=.75, fill=black, label=left:]{};
    \vertex (03) at (1, 7.5)  [scale=.75, fill=black, label=right:]{};
    \vertex (04) at (1, 9.5)  [scale=.75, fill=black, label=above:]{};
    \vertex (05) at (1.5, 8.5)  [scale=.75, fill=black, label=right:]{};
    \vertex (06) at (1.5, 6.5)  [scale=.75, fill=black, label=below:]{};
    \vertex (07) at (.5, 6.5)  [scale=.75, fill=black, label=below:]{};
     \vertex (08) at (2, 7)  [scale=.75, fill=black, label=below:]{};
  
     \vertex (0B1) at (4.5,7)  [scale=.75, fill=black, label=below:]{};
    \vertex (0B2) at (4.5,8.5)  [scale=.75, fill=black, label=left:]{};
    \vertex (0B3) at (5, 7.5)  [scale=.75, fill=black, label=below:]{};
    \vertex (0B4) at (5, 9.5)  [scale=.75, fill=black, label=above:]{};
    \vertex (0B5) at (5.5, 8.5)  [scale=.75, fill=black, label=right:]{};
    \vertex (0B6) at (6.5, 6)  [scale=.75, fill=black, label=below:]{};
    \vertex (0B7) at (6, 6.5)  [scale=.75, fill=black, label=below:]{};
    \vertex (0B8) at (5.5, 7)  [scale=.75, fill=black, label=below:]{};

       \vertex (0C1) at (9,7)  [scale=.75, fill=black, label=below:]{};
    \vertex (0C2) at (9,8.5)  [scale=.75, fill=black, label=left:]{};
    \vertex (0C3) at (9.5, 7.5)  [scale=.75, fill=black, label=below:]{};
    \vertex (0C4) at (9.5, 9.5)  [scale=.75, fill=black, label=above:]{};
    \vertex (0C5) at (10, 8.5)  [scale=.75, fill=black, label=right:]{};
    \vertex (0C6) at (10, 6)  [scale=.75, fill=black, label=below:]{};
    \vertex (0C7) at (10.75, 6.5)  [scale=.75, fill=black, label=right:]{};
    \vertex (0C8) at (10, 7)  [scale=.75, fill=black, label=right:]{};
    
      \vertex (1) at (13.5,9.5)  [scale=.75, fill=black, label=above:]{};
    \vertex (2) at (14.5, 9.5)  [scale=.75, fill=black, label=above:]{};
    \vertex (3) at (14, 8.5)  [scale=.75, fill=black, label=right:]{};
    \vertex (4) at (14, 7.5)  [scale=.75, fill=black, label=right:]{};
    \vertex (5) at (13.25, 6.5)  [scale=.75, fill=black, label=below:]{};
    \vertex (6) at (13.75, 6.5)  [scale=.75, fill=black, label=below:]{};
    \vertex (7) at (14.25, 6.5)  [scale=.75, fill=black, label=below:]{};
     \vertex (8) at (14.75, 6.5)  [scale=.75, fill=black, label=below:]{};

    \node(A) at (1, 4.75) []{$V_1$};
    \node(B) at (5, 4.75)[]{$V_2$};
    \node(C) at (9.5, 4.75)[]{$V_3$};
    \node(A) at (14, 4.75) []{$V_4$};

     \vertex (B1) at (.5,2.5)  [scale=.75, fill=black, label=below:]{};
    \vertex (B2) at (1.5,2.5)  [scale=.75, fill=black, label=left:]{};
    \vertex (B3) at (1, 1.5)  [scale=.75, fill=black]{};
    \vertex (B4) at (.5, .5)  [scale=.75, fill=black, label=above:]{};
    \vertex (B5) at (1.5, .5)  [scale=.75, fill=black, label=right:]{};
    \vertex (B6) at (.5, -.5)  [scale=.75, fill=black, label=below:]{};
    \vertex (B7) at (1.25, -.5)  [scale=.75, fill=black, label=below:]{};
       \vertex (B8) at (1.75, -.5)  [scale=.75, fill=black, label=below:]{};

        \vertex (C1) at (4.5,2.5)  [scale=.75, fill=black, label=below:]{};
    \vertex (C2) at (5.5,2.5)  [scale=.75, fill=black, label=left:]{};
    \vertex (C3) at (5, 1.5)  [scale=.75, fill=black]{};
    \vertex (C4) at (5, 1)  [scale=.75, fill=black, label=above:]{};
    \vertex (C5) at (4.5, .5)  [scale=.75, fill=black, label=right:]{};
    \vertex (C6) at (5.5, .5)  [scale=.75, fill=black, label=below:]{};
    \vertex (C7) at (5.5, 0)  [scale=.75, fill=black, label=below:]{};
        \vertex (C8) at (5.5, -.5)  [scale=.75, fill=black, label=below:]{};

        \vertex (D1) at (9,2.5)  [scale=.75, fill=black, label=below:]{};
    \vertex (D2) at (10,2.5)  [scale=.75, fill=black, label=left:]{};
    \vertex (D3) at (9, 1.5)  [scale=.75, fill=black]{};
    \vertex (D4) at (10, 1.5)  [scale=.75, fill=black, label=above:]{};
    \vertex (D5) at (9.5, .5)  [scale=.75, fill=black, label=right:]{};
    \vertex (D8) at (9.5, 0)  [scale=.75, fill=black, label=below:]{};
    \vertex (D6) at (9, -.5)  [scale=.75, fill=black, label=below:]{};
    \vertex (D7) at (10, -.5)  [scale=.75, fill=black, label=below:]{};
    
     \vertex (E1) at (13.5,2.5)  [scale=.75, fill=black, label=below:]{};
    \vertex (E2) at (14.5,2.5)  [scale=.75, fill=black, label=left:]{};
    \vertex (E3) at (13.5, 1.5)  [scale=.75, fill=black]{};
    \vertex (E4) at (14.5, 1.5)  [scale=.75, fill=black, label=above:]{};
    \vertex (E5) at (14, .5)  [scale=.75, fill=black, label=right:]{};
    \vertex (E6) at (13.5, -.5)  [scale=.75, fill=black, label=below:]{};
    \vertex (E7) at (14.5, -.5)  [scale=.75, fill=black, label=below:]{};
      \vertex (E8) at (15, .5)  [scale=.75, fill=black, label=below:]{};
    
         \vertex (F1) at (.5,-4)  [scale=.75, fill=black, label=below:]{};
    \vertex (F2) at (1.5,-4)  [scale=.75, fill=black, label=left:]{};
    \vertex (F3) at (.5, -5)  [scale=.75, fill=black]{};
    \vertex (F4) at (1.5, -5)  [scale=.75, fill=black, label=above:]{};
    \vertex (F5) at (.5, -6)  [scale=.75, fill=black, label=right:]{};
    \vertex (F6) at (1.5, -6)  [scale=.75, fill=black, label=below:]{};
    \vertex (F7) at (1, -7)  [scale=.75, fill=black, label=below:]{};
      \vertex (F8) at (-.5, -5)  [scale=.75, fill=black, label=below:]{};

         \vertex (G1) at (4.5,-4)  [scale=.75, fill=black, label=below:]{};
    \vertex (G2) at (5.5,-4)  [scale=.75, fill=black, label=left:]{};
    \vertex (G3) at (5, -5)  [scale=.75, fill=black]{};
    \vertex (G4) at (4.5, -6)  [scale=.75, fill=black, label=above:]{};
    \vertex (G5) at (5.5, -6)  [scale=.75, fill=black, label=right:]{};
    \vertex (G6) at (5, -6.5)  [scale=.75, fill=black, label=below:]{};
    \vertex (G7) at (4.75, -7)  [scale=.75, fill=black, label=below:]{};
        \vertex (G8) at (5.25, -7)  [scale=.75, fill=black, label=below:]{};

         \vertex (H1) at (9,-4)  [scale=.75, fill=black, label=below:]{};
    \vertex (H2) at (10,-5)  [scale=.75, fill=black, label=left:]{};
    \vertex (H3) at (9, -5)  [scale=.75, fill=black]{};
    \vertex (H4) at (10, -6)  [scale=.75, fill=black, label=above:]{};
    \vertex (H5) at (9, -6)  [scale=.75, fill=black, label=right:]{};
    \vertex (H6) at (10, -7)  [scale=.75, fill=black, label=below:]{};
    \vertex (H7) at (9, -7)  [scale=.75, fill=black, label=below:]{};
       \vertex (H8) at (11, -6)  [scale=.75, fill=black, label=below:]{};
    
      \vertex (I1) at (13.5,-4)  [scale=.75, fill=black, label=below:]{};
    \vertex (I2) at (14.5,-4)  [scale=.75, fill=black, label=left:]{};
    \vertex (I3) at (14, -5)  [scale=.75, fill=black]{};
    \vertex (I4) at (14, -5.5)  [scale=.75, fill=black, label=above:]{};
    \vertex (I5) at (14, -6)  [scale=.75, fill=black, label=right:]{};
    \vertex (I6) at (13.5, -7)  [scale=.75, fill=black, label=below:]{};
    \vertex (I7) at (14.5, -7)  [scale=.75, fill=black, label=below:]{};
        \vertex (I8) at (15, -5.5)  [scale=.75, fill=black, label=below:]{};

    \node(B) at (1, -2)[]{$V_5$};
    \node(C) at (5, -2)[]{$V_6$};
     \node(D) at (9.5, -2)[]{$V_7$};
      \node(E) at (14, -2)[]{$V_8$};
       \node(F) at (1, -8.5)[]{$V_9$};
         \node(G) at (5, -8.5)[]{$V_{10}$};
           \node(H) at (9.5, -8.5)[]{$V_{11}$};
             \node(I) at (14, -8.5)[]{$V_{12}$};

    \path 
    
    	(01) edge (03)
	(06) edge (03)
	(07) edge (03)
	(08) edge (03)
	(02) edge (03)
	(04) edge (05)
	(02) edge (04)
	(03) edge (05)
	(02) edge(05)
	
	(0B1) edge (0B3)
	(0B2) edge (0B3)
	(0B2) edge (0B5)
	(0B4) edge (0B5)
	(0B3) edge (0B5)
	(0B4) edge (0B2)
	(0B3) edge (0B8)
	(0B8) edge (0B7)
	(0B6) edge (0B7)

  	(0C1) edge (0C3)
	(0C2) edge (0C3)
	(0C2) edge (0C5)
	(0C4) edge (0C5)
	(0C3) edge (0C5)
	(0C4) edge (0C2)
	(0C3) edge (0C8)
	(0C8) edge (0C7)
	(0C6) edge (0C7)
	(0C8) edge (0C6)
	(1) edge (2)
	(1) edge (3)
	(2) edge (3)
	(4) edge (3)
	(4) edge (5)
	(4) edge (6)
	(4) edge (7)
	(4) edge (8)
	
	(B1) edge (B2)
	(B1) edge (B3)
	(B2) edge (B3)
	(B3) edge (B4)
	(B3) edge (B5)
	(B4) edge (B6)
	(B5) edge (B7)
	(B5) edge (B8)

	(C1) edge (C2)
  	(C1) edge (C3)
	(C2) edge (C3)
	(C3) edge (C4)
	(C4) edge (C5)
	(C4) edge (C6)
	(C6) edge (C7)
	(C7) edge (C8)

		(D1) edge (D2)
  	(D1) edge (D3)
	(D2) edge (D4)
	(D3) edge (D4)
	(D4) edge (D5)
	(D3) edge (D5)
	(D5) edge (D8)
	(D8) edge (D7)
	(D8) edge (D6)

		(E1) edge (E2)
  	(E1) edge (E3)
	(E2) edge (E4)
	(E3) edge (E4)
	(E4) edge (E5)
	(E4) edge (E8)
	(E5) edge (E6)
	(E5) edge (E7)
	(E6) edge (E7)
	
		(F1) edge (F2)
  	(F1) edge (F3)
	(F2) edge (F4)
	(F3) edge (F4)
	(F4) edge (F6)
		(F1) edge (F8)
	(F5) edge (F6)
	(F5) edge (F7)
	(F6) edge (F7)
	
	(G1) edge (G2)
  	(G1) edge (G3)
	(G2) edge (G3)
	(G3) edge (G4)
	(G3) edge (G5)
		(G4) edge (G5)
	(G3) edge (G6)
	(G6) edge (G7)
	(G6) edge (G8)

(H3) edge (H2)
(H4) edge (H2)
	(H5) edge (H3)
  	(H5) edge (H4)
	(H5) edge (H6)
	(H7) edge (H5)
	(H6) edge (H7)
	(H3) edge (H1)
	(H4) edge (H8)
	
	(I1) edge (I2)
  	(I1) edge (I3)
	(I2) edge (I3)
	
  	(I5) edge (I6)
	(I5) edge (I7)
	(I6) edge (I7)
	(I4) edge (I8)
	(I3) edge (I4)
	(I5) edge (I4)
	
    ;

\end{tikzpicture}
\caption{All connected, $K_4$-free, well-edge dominated graphs containing a triangle of order $8$ where only $V_1, V_2$, and $V_3$  contain a diamond}
\label{ND8}
\end{center}
\end{figure}
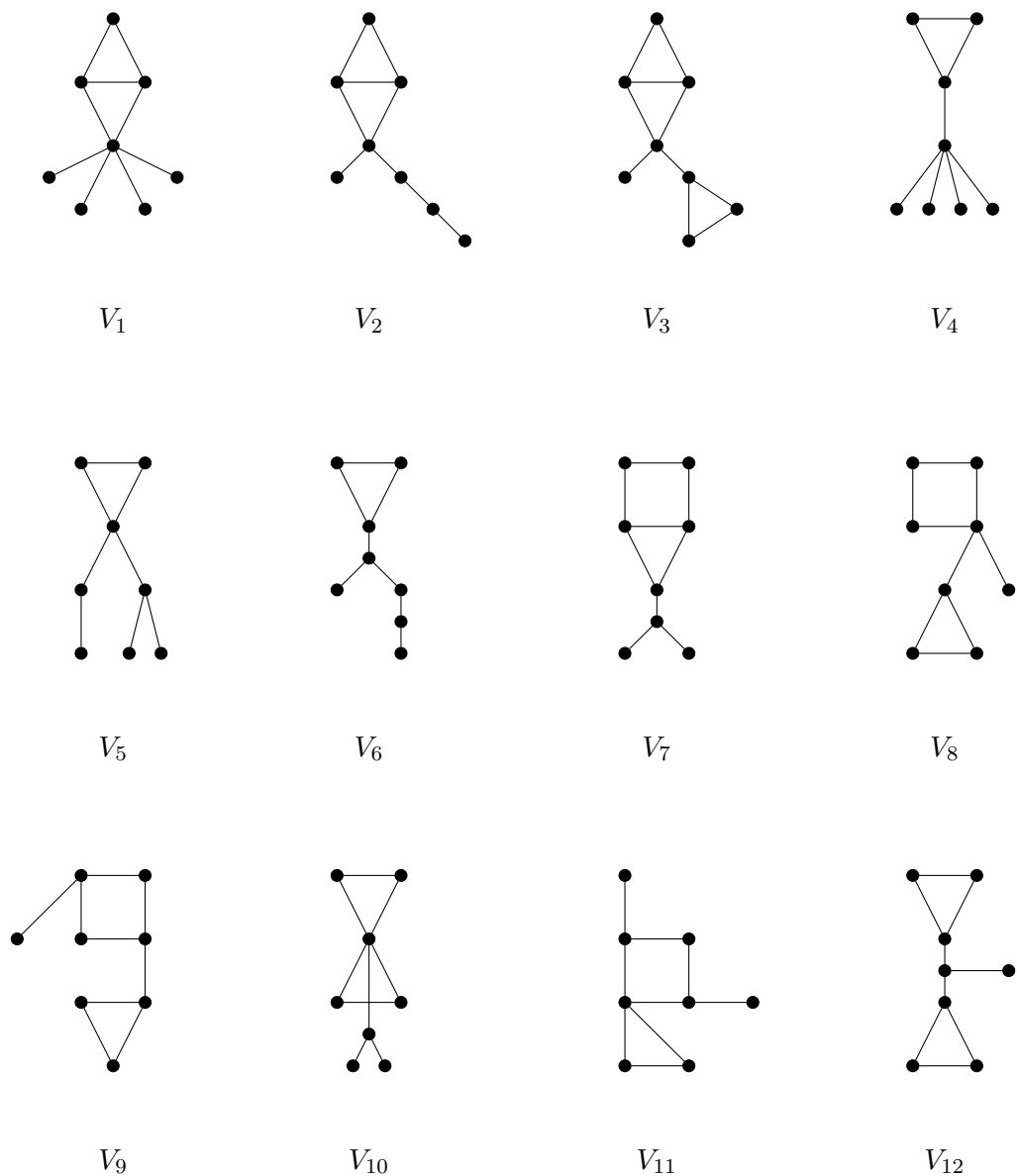

Our first step in finding all $K_4$-free well-edge-dominated graphs was to have Sage find all such graphs up to order $8$. Figures~\ref{ND7} and \ref{ND8} are a complete list of all connected, $K_4$-free, nonbipartite,  well-edge-dominated graphs containing a triangle of orders $7$ and $8$, respectively. 
We are now ready to define three infinite classes of graphs, each of which whose members are well-edge-dominated. First, we define the class $\mathcal{P}$ (referred to as the propellers) and the class $\mathcal{W}$ (referred to as the windmills). Let $T_1, \dots, T_k$ be $k$ disjoint triangles and for each $i \in [k]$ label one vertex of triangle $T_i$ by $v_i$. $G \in \mathcal{P}$ if it is obtained from $T_1, \dots, T_k$ by identifying each vertex in $\{v_1, \dots, v_k\}$ into a single vertex (see $W_{13}$ in Figure~\ref{ND7}). To define $\mathcal{W}$, we let $\mathcal{H}$ be the house graph, depicted in Figure~\ref{fig:houses}(b), and label the vertex on the triangle with degree two $x$. $G \in \mathcal{W}$ if it is obtained from $T_1, \dots, T_k,$ and $ \mathcal{H}$ by identifying each vertex in $\{v_1, \dots, v_k, x\}$ into a single vertex (see $W_8$ in Figure~\ref{ND7}). In either case, the identified vertex in $G$ will be referred to as the ``nose" of the windmill/propeller.  We will also consider  $K_3$ as a propeller and $\mathcal{H}$ as a windmill. Before defining our broader class, we show that every graph in $\mathcal{P}\cup \mathcal{W}$ is indeed well-edge-dominated. 

\begin{lemma}\label{lem:windmill} If $G \in \mathcal{P}\cup\mathcal{W}$, then $G$ is well-edge-dominated. 
\end{lemma}

\begin{proof} Suppose first that $G \in \mathcal{P}$ and is obtained from the triangles $T_1, \dots, T_k$. Any minimal edge dominating set will contain exactly one edge from each of the $k$ edge disjoint triangles in $G$, meaning that $\gamma'(G) = k = \Gamma'(G)$. 

Next, let $G \in \mathcal{W}$. Note that if $G \cong \mathcal{H}$, then $G$ is well-edge-dominated. Thus, we may assume $G \not\cong \mathcal{H}$. In addition, assume we have obtained $G$ from the triangles $T_1, \dots, T_k$ and the house graph $\mathcal{H}$ where the vertex in $\mathcal{H}$ on the triangle with degree two is $x$ by identifying each vertex in $\{v_1, \dots, v_k, x\}$ (for any $v_i$ in $T_i$, $i\in[k]$) into a single vertex $w$. Any minimal edge dominating set $D$ of $G$ will contain exactly one edge from each of the $k$ triangles in $G$ originally found in $T_1, \dots, T_k$. Moreover, no matter which $k$ edges are chosen from $T_1, \dots, T_k$, the vertices in $\mathcal{H} - x$ induce a $C_4$. So $D$ must contain at least two edges from $\mathcal{H}$. To see that $D$ contains exactly two edges from $\mathcal{H}$, suppose to the contrary that $D$ contains three edges from $\mathcal{H}$. Then all three edges must be incident to a vertex with degree three in $\mathcal{H}$, say $y$. However, in this case, $xy$ must have a private edge neighbor which is not in $\mathcal{H}$ and since all edges outside $\mathcal{H}$ are dominated by an edge in its respective triangle, $D$ is not minimal. Therefore, $D$ contains exactly two edges from $\mathcal{H}$. It follows that $G$ is well-edge-dominated. 
\end{proof}

Now, we are ready to define our broader class $\mathcal{G}$, built from the classes $\mathcal{P} \cup \mathcal{W}$, as follows. Consider the graph obtained from the disjoint union of a diamond and $K_2$ by identifying a vertex of $K_2$ with an exterior vertex of the diamond. Call the graph $K^*$ and put $K^*$ in $\mathcal{G}$. Additionally, if $G \in \mathcal{P}\cup \mathcal{W}$, then put $G$ in $\mathcal{G}$. All other graphs $G \in \mathcal{G}$ can be constructed in the following way. Let $G'=(A \cup B, E')$ be a connected, nontrivial, bipartite and well-edge-dominated graph with $|A| < |B|$.  Choose a set $B' \subset B$ and a set $Y \subseteq A$ such that the following are true:
\begin{itemize}
\item $B'\cup Y \ne \emptyset$.
\item  $G'- B'$ is a well-edge-dominated graph with $\gamma'(G'-B') = \gamma'(G')$, $G'-B'$ has no trivial components, and $|A| \le |B-B'|$. 
\item If $Y\ne \emptyset$, then we write $Y = \{y_1, \dots, y_{\ell}\}$ where $y_i$ is a support vertex in $G'- B'$. 
\end{itemize}
 Each vertex $v$ in $B'$ is called detachable and it is called strongly detachable if each neighbor of it in $G'$ is a support vertex in $G'-B'$.  If $Y\ne \emptyset$, let $D_1, \dots, D_{\ell}$ be a disjoint union of diamonds, and for each $D_i$, label one of the vertices of degree two as the ``nose" of the diamond. If $B' \ne\emptyset$, then choose a set $B'' \subseteq B'$ (possibly empty) where each vertex in $B''$ is strongly detachable. If $B'' \ne \emptyset$, then write $B'' = \{s_1, \dots, s_k\}$ and choose any set of $k$ windmills, enumerated as $W_1, \dots, W_k$. If $B' - B'' \ne \emptyset$, write $B' - B'' = \{x_1, \dots, x_r\}$ and choose any set of $r$ propellers, enumerated as $P_1, \dots, P_r$.  We obtain $G$ from $G'$ under the following rules: 
\begin{enumerate}
\item[(a)] The nose of $W_i$ is identified with $s_i$ in $B''$.
\item[(b)] The nose of $P_i$ is identified with $x_i$ in $B' - B''$. 
\item[(c)] The nose of $D_i$ is identified with the support vertex $y_i$ of $G'$ (which stays a support vertex in $G$). 
\end{enumerate}

Note that of all the graphs shown in Figures~\ref{ND7} and \ref{ND8}, only $W_1, W_2, W_3, W_{10}$, and $W_{12}$ are not in $\mathcal{G}$. However, $W_{10} \cong\mathcal{DH}$ and $W_{12} \cong Cr$. So there are only three exceptions of connected, well-edge-dominated graphs with order $7$ or $8$ that both contain at least two triangles (and are therefore not in $\mathcal{T}\cup \mathcal{F}$) and are not in $\mathcal{G}$. We next prove that each graph in $\mathcal{G}$ is in fact well-edge-dominated. First, we prove a more general result about constructing well-edge-dominated graphs from two smaller well-edge-dominated graphs. 

\begin{lemma}\label{lem:construct} Let $G_1$ and $G_2$ be two nontrivial well-edge-dominated graphs with $x \in V(G_1)$ and $y \in V(G_2)$ such that $G_1-x$ is well-edge-dominated with $\gamma'(G_1-x) = \gamma'(G_1)$ and $G_2-y$ is well-edge-dominated with $\gamma'(G_2-y) = \gamma'(G_2)$. Then the graph $H$ obtained from $G_1$ and $G_2$ by identifying $x$ and $y$ is well-edge-dominated with $\gamma'(H) = \gamma'(G_1) + \gamma'(G_2)$. 
\end{lemma}
\begin{proof} Let $w$ be the identified vertex in $H$ and let $F$ be a minimal edge dominating set of $H$. We claim that $F \cap E(G_1 - x)$ is a minimal edge dominating set of $G_1 - x$. Indeed, if not, then some edge $e \in F \cap E(G_1-x)$ has no private edge neighbor in $G_1 - x$ with respect to $F \cap E(G_1 - x)$. However, $e$ does have a private edge neighbor in $H$ with respect to $F$, meaning $e$ has a private edge neighbor incident to $w$. This in turn implies that no edge in $F$ is incident to $w$ and $F \cap E(G_1-x)$ is a minimal edge dominating set of $G_1$. Hence, $|F \cap E(G_1- x)| = \gamma'(G_1) = \gamma'(G_1-x)$, a contradiction. It follows that $F \cap E(G_1-x)$ is a minimal edge dominating set of $G_1 - x$, and as $G_1 - x$ is well-edge-dominated, $|F \cap E(G_1-x)| = \gamma'(G_1-x) = \gamma'(G_1)$. Similarly, $|F \cap E(G_2 - y)| = \gamma'(G_2-y) = \gamma'(G_2)$ from which it follows that $|F| = \gamma'(G_1) + \gamma'(G_2)$. As $F$ is an arbitrary minimal edge dominating set of $H$, $H$ is well-edge-dominated. 
\end{proof}

Note that each graph in $\mathcal{G}$ can be viewed as performing a sequence of identifications between two well-edge-dominated graphs. Next,  we show every graph in $\mathcal{G}$ is in fact well-edge-dominated. 

\begin{theorem}\label{cor:G} If $G \in \mathcal{G}$, then $G$ is well-edge-dominated. 
\end{theorem}

\begin{proof}
Note that the graph obtained by identifying a vertex of $K_2$ with an exterior vertex of a diamond is well-edge-dominated. Moreover, if $G \in \mathcal{P}\cup \mathcal{W}$, then Lemma~\ref{lem:windmill} shows $G$ is well-edge-dominated. Therefore, we shall assume that $G$ was constructed from the connected, nontrivial, bipartite and well-edge-dominated graph $G'= (A\cup B, E')$ where $|A| < |B|$. Then a set $B' \subset B$ and $Y \subseteq A$ was chosen so that  \begin{itemize}
\item $B'\cup Y \ne \emptyset$.
\item  $G'- B'$ is a well-edge-dominated graph with $\gamma'(G'-B') = \gamma'(G')$, $G'-B'$ has no trivial components, and $|A| \le |B-B'|$. 
\item If $Y\ne \emptyset$, then we write $Y = \{y_1, \dots, y_{\ell}\}$ where $y_i$ is a support vertex in $G'- B'$. 
\end{itemize}
If $Y\ne \emptyset$, we enumerate $Y = \{y_1, \dots, y_{\ell}\}$ and let $D_1, \dots, D_{\ell}$ be the disjoint set of diamonds such that the nose of $D_i$ was identified with $y_i$. If $B'' \ne \emptyset$, we enumerate $B'' = \{s_1, \dots, s_k\}$ and let $W_1, \dots, W_k$ be the disjoint set of windmills such that the nose of $W_i$ was identified with $s_i$. If $B' - B'' \ne \emptyset$, we enumerate $B' - B'' = \{x_1, \dots, x_r\}$ and let $P_1, \dots, P_r$ be the set of disjoint propellers such that the nose of $P_i$ was identified with $x_i$.  In $G$, we give the vertex representing an identification with a vertex in $B' \cup Y$ the same label as the vertex in $B' \cup Y$. 

Suppose first that $W_i$ was used to create $G$ where $W_i$ is a windmill that was created from $n_i$ triangles and the house graph. We claim that $|E(W_i)\cap F| = n_i + 2$. Note first that each of the $n_i$ triangles must contain an edge from $F$. Moreover, the house graph contains a $4$-cycle whose vertices are not the nose of $W_i$. This means that $F$ contains at least two edges from the house graph in $W_i$. To see that $F$ does not contain three edges from the house, label the vertices on the house graph in $W_i$ as $\{x_i, y_i, s_i, c_i, d_i\}$ where $s_i$ is the nose of $W_i$, $x_iy_is_i$ is the triangle, and $c_ix_i$ and $d_iy_i$ are edges in $W_i$. Note that the only way for $F$ to contain $3$ edges from the house graph is if $F$ contains $\{x_is_i, x_iy_i, a_ix_i\}$ or $F$ contains $\{y_is_i, y_ix_i, y_ib_i\}$. Without loss of generality, we may assume $F$ contains $\{x_is_i, x_iy_i, a_ix_i\}$. It follows that the private edge neighbor of $x_is_i$ must be an edge incident to $s_i$, but not $y_is_i$. Moreover, since each of the $n_i$ triangles used to create $W_i$ contain an edge of $F$, the private edge neighbor of $x_is_i$ is $s_iw_i$ where $s_iw_i$ was originally an edge in $G'$. Moreover, by the construction of $G$, $s_i$ is a strongly detachable vertex in $B'$, meaning that  $s_i$  is a support vertex in $G' - B'$ so that $F$ contains an edge incident to $s_i$. Therefore, $x_is_i$ indeed has no private edge neighbor and $|E(W_i) \cap F| = n_i+2$.

Next, suppose for some $i \in [r]$ that $P_i$ was used to create $G$ where $P_i$ is a propeller that was created from $n_i$ triangles. It is clear that $|E(P_i) \cap F| = n_i$. Lastly, suppose for some $i \in [\ell]$ that $D_i$ was used to create $G$ where $D_i$ is a diamond with vertices $x_i, y_i, z_i, w_i$ where $y_i$ is the nose of the diamond and $w_i$ is the other exterior vertex of  the diamond. Note that $|E(D_i) \cap F| \le 2$. Suppose that $|E(D_i)| \cap F| = 1$. It follows that $x_iz_i \in F$ for otherwise all edges of the triangle $w_ix_iz_i$ are not edge dominated. Note further that by construction of $G$, $y_i$ is a support vertex in $G'-B'$ so that $F$ contains an edge of the form $y_ia$ where $y_ia$ is an edge in $G'$. We may weakly partition the diamonds $D_1, \dots, D_{\ell}$ as $D_1, \dots, D_{\alpha}, D_{\alpha +1}, \dots, D_{\ell}$ where $|E(D_i) \cap F| = 2$ for $i \in [\alpha]$ and $|E(D_i)\cap F| = 1$ otherwise. 

Next, note that if $|E(D_i)\cap F|=2$, then at least one of the edges in $E(D_i)\cap F$ is incident to $y_i$. This is because in $G'-B'$, $y_i$ is a support vertex. Therefore, if neither edge in $E(D_i) \cap F$ is incident to $y_i$, then $E(G') \cap F$ contains an edge incident to $y_i$ which would imply that for some edge in $E(D_i) \cap F$, say $e$, $F- \{e\}$ is also an edge dominating set of $G$, contradicting the minimality of $F$. Now for each $i \in [\alpha]$, let $a_i$ be a leaf of $y_i$ in $G'-B'$, $Z$ be the set of all edges of the form $a_iy_i$ for $i \in [\alpha]$, and let $F'\subseteq F$ be a matching that saturates all $y_i \in Y$ where $\alpha+1 \le i \le \ell$. As $G'$ is well-edge-dominated, $G' - N_e[F'\cup Z]$ is well-edge-dominated by Lemma~\ref{lem:reduce}. Thus, any minimal edge dominating set of $G'-N_e[F'\cup Z]$ has cardinality $|A| - \ell$. Therefore, there are exactly $|A|-\ell$ edges in $E(G') \cap F$ that are not incident to a vertex in $Y$. It follows that 
\[|E(G') \cap F| + \sum_{i=1}^{\ell} |E(D_i)\cap F| \ge |A|- \ell + 2\ell = |A|+\ell.\]
Seeing as how $F$ is arbitrary, we have that 
\[|F| = |E(G')\cap F| + \sum_{i=1}^k (n_i +2) + \sum_{i=1}^r n_i + \sum_{i=1}^{\ell} |E(D_i)\cap F| = |A| +\ell + \sum_{i=1}^k (n_i +2) + \sum_{i=1}^r n_i \]
and $G$ is well-edge-dominated.

\end{proof}

The remainder of the paper is now focused on showing that we have found all $K_4$-free well-edge-dominated graphs of girth $3$. We note that our proof will rely heavily on considering when a well-edge-dominated graph contains a diamond or not. For this reason, we point out the following observation about the graphs depicted in Figures~\ref{ND7} and \ref{ND8}. 
\begin{obs}
\label{obs:diamondsize7}
 Assume $G$ has order $7$ and contains a diamond as an induced subgraph. $G$ is a $K_4$-free well-edge-dominated graph with girth $3$ if and only if $G \in \{W_1, W_2, W_3, W_4\}$, where $W_4$ is the graph in $\mathcal{G}$ that contains a diamond as well as three pendent edges.
\end{obs}

\begin{figure}[h!]
\begin{center}
\begin{tikzpicture}[scale=.75]
    \vertex (1) at (0,0)  [scale=.75, fill=black, label=below:$a$]{};
    \vertex (2) at (0,1)  [scale=.75, fill=black, label=left:$w$]{};
    \vertex (3) at (1, .5)  [scale=.75, fill=black, label=below:$y$]{};
    \vertex (4) at (1, 1.5)  [scale=.75, fill=black, label=above:$x$]{};
    \vertex (5) at (2, 1)  [scale=.75, fill=black, label=right:$z$]{};
    \vertex (6) at (2, 0)  [scale=.75, fill=black, label=below:$c$]{};
    \vertex (7) at (1, -.5)  [scale=.75, fill=black, label=below:$b$]{};
  
     \vertex (B1) at (4,0)  [scale=.75, fill=black, label=below:$a$]{};
    \vertex (B2) at (4,1)  [scale=.75, fill=black, label=left:$w$]{};
    \vertex (B3) at (5, .5)  [scale=.75, fill=black, label=below:$y$]{};
    \vertex (B4) at (5, 1.5)  [scale=.75, fill=black, label=above:$x$]{};
    \vertex (B5) at (6, 1)  [scale=.75, fill=black, label=right:$z$]{};
    \vertex (B6) at (6, 0)  [scale=.75, fill=black, label=below:$c$]{};
    \vertex (B7) at (5, -.5)  [scale=.75, fill=black, label=below:$b$]{};
    
       \vertex (C1) at (8,0)  [scale=.75, fill=black, label=below:$a$]{};
    \vertex (C2) at (8,1)  [scale=.75, fill=black, label=left:$w$]{};
    \vertex (C3) at (9, .5)  [scale=.75, fill=black, label=below:$y$]{};
    \vertex (C4) at (9, 1.5)  [scale=.75, fill=black, label=above:$x$]{};
    \vertex (C5) at (10, 1)  [scale=.75, fill=black, label=right:$z$]{};
    \vertex (C6) at (10, 0)  [scale=.75, fill=black, label=below:$c$]{};
    \vertex (C7) at (9, -.5)  [scale=.75, fill=black, label=below:$b$]{};

       \vertex (D1) at (12,0)  [scale=.75, fill=black, label=below:$a$]{};
    \vertex (D2) at (12,1)  [scale=.75, fill=black, label=left:$w$]{};
    \vertex (D3) at (13, .5)  [scale=.75, fill=black, label=below:$y$]{};
    \vertex (D4) at (13, 1.5)  [scale=.75, fill=black, label=above:$x$]{};
    \vertex (D5) at (14, 1)  [scale=.75, fill=black, label=right:$z$]{};
    \vertex (D6) at (14, 0)  [scale=.75, fill=black, label=below:$c$]{};
    \vertex (D7) at (13, -.5)  [scale=.75, fill=black, label=below:$b$]{};

    \node(A) at (1, -1.75) []{$W_1$};
    \node(B) at (5, -1.75)[]{$W_2$};
    \node(C) at (9, -1.75)[]{$W_3$};
     \node(D) at (13, -1.75)[]{$W_4$};

    \path 
	(1) edge (2)
	(2) edge (3)
	(3) edge (4)
	(4) edge (5)
	(5) edge (6)
	(6) edge (7)
	(1) edge (7)
	(2) edge (4)
	(3) edge (5)
	
	(B1) edge (B2)
	(B2) edge (B3)
	(B3) edge (B4)
	(B4) edge (B5)
	(B5) edge (B6)
	(B6) edge (B7)
	(B1) edge (B7)
	(B2) edge (B4)
	(B3) edge (B5)
	(B3) edge [bend left]  (B7)
	
  	(C1) edge (C2)
	(C2) edge (C3)
	(C3) edge (C4)
	(C4) edge (C5)
	(C5) edge (C6)
	(C6) edge (C7)
	(C1) edge (C7)
	(C2) edge (C4)
	(C3) edge (C5)  
	(C4) edge (C6)
	
	(D1) edge (D2)
	(D2) edge [bend left] (D6)
	(D2) edge  (D7)
	(D2) edge (D3)
	(D3) edge (D4)
	(D4) edge (D5)
	(D2) edge (D4)
	(D3) edge (D5)  

    ;

\end{tikzpicture}
\caption{The graphs $W_1$, $W_2$, $W_3$, and $W_4$}
\label{fig:Ws}
\end{center}
\end{figure}
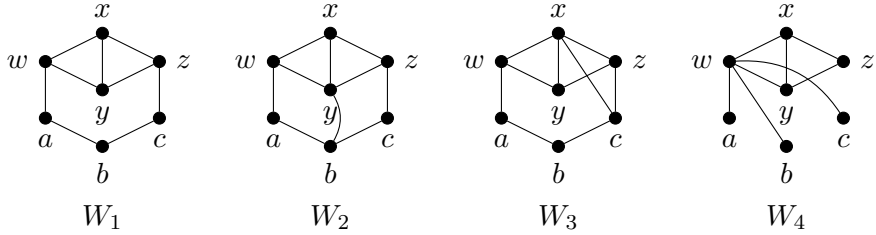

%\begin{obs}
%\label{obs:diamondsize8}
%Let $G$ be $K_4$-free. Assume $G$ has order $8$ and contains a diamond as an induced subgraph. $G$ is well-edge-dominated with girth $3$ if and only if $G \in \{V_1, V_2, V_3\}$, which are shown in Figure \ref{fig:Vs}.
%\end{obs}

Before we prove the general result, we consider well-edge-dominated graphs containing triangles that have small diameter. We refer to any graph $G$ obtained by appending $n\ge 1$ leaves to exactly one exterior vertex of a diamond as a \emph{kite}. We let $\mathcal{K}$ be the set of all kites. We also let $\mathcal{R}$ be the set of all graphs $G$ obtained from the disjoint union of a propeller or windmill and $k\ge 1$ copies of $P_3$ by identifying the nose of the propeller or windmill with exactly one leaf from each of the $P_3$s. Note that $\mathcal{K}\cup \mathcal{R} \subset \mathcal{G}$.

\begin{lemma}\label{lem:small} Let $G$ be a connected $K_4$-free well-edge-dominated graph built from an independent set $I$, a set $S$, and a triangle $xyz$ such that  no vertex $t \in I$ is adjacent to a vertex in $\{x, y, z\}\cup (I - \{t\})$, and each vertex in $S$ has a neighbor in $\{x, y, z\}$. Then $G \in \{Cr, \mathcal{H}\}\cup \mathcal{K} \cup \mathcal{P} \cup \mathcal{R}\cup \mathcal{W}$.
\end{lemma}

\begin{proof}
Suppose to the contrary that a counterexample exists and among all such counterexamples, choose one with minimum cardinality. We may assume that $n(G) \ge 9$. Furthermore, if $G$ contains only one triangle, then by Theorem~\ref{thm:onetriangle} and the assumptions on $G$, $G \in \mathcal{F} \cup \mathcal{T}\cup \{\mathcal{H}, K_3, Cr\}$ where $K_3 \in \mathcal{P}$ and $\mathcal{H} \in \mathcal{W}$. Thus, we may assume that $G$ contains at least two triangles. We proceed by considering whether $I= \emptyset$ or $I \ne \emptyset$. 
\vskip2mm
\noindent\textbf{Case 1:} Suppose $I\ne \emptyset$.
\vskip2mm
Let $st \in E(G)$ where $t \in I$ and $s \in S$. As $G$ is assumed to be well-edge-dominated, $H = G - N_e[st]$ is well-edge-dominated. Moreover, $H$ satisfies the assumptions in the statement of the theorem so we may conclude that $H \in \{Cr, \mathcal{H}\}\cup \mathcal{K}\cup \mathcal{P} \cup \mathcal{W}\cup \mathcal{R}$.

Assume first that at least two vertices of $\{x, y, z\}$ have degree at least three in $H$. Without loss of generality, we may assume $\deg_H(x) \ge 3$ and $\deg_H(y) \ge 3$. It follows that $H \in \{Cr, \mathcal{H}\} \cup \mathcal{K}\cup \mathcal{W}\cup  \mathcal{R}$. However, we have assumed $n(H) \ge 7$ meaning that $H \cong Cr$ or $H \in \mathcal{K}\cup \mathcal{W}\cup \mathcal{R}$. Let us assume that $H \cong Cr$ and label the vertices as in Figure~\ref{fig:houses} (c). It follows that $\{s, a, b, c\} = S$ and $\{t, d\} = I$. Without loss of generality, we may assume that $s$ is adjacent $x$ or $y$. If $sx\in E(G)$, then  $H' = G - N_e[cd]$ is a well-edge-dominated graph containing a triangle with order $7$ and yet by choice of minimal counterexample, $H' \in \{Cr, \mathcal{H}\} \cup \mathcal{K} \cup \mathcal{P}\cup \mathcal{R}\cup \mathcal{W}$. However, this cannot be as $\deg_{H'}(x) \ge 4$, $\deg_{H'}(y) \ge 3$, and $\{x, y, z, a, b\}$ induces the house graph and no such graph in $\mathcal{K}\cup \mathcal{P}\cup \mathcal{R}\cup \mathcal{W}$ satisfies these conditions. If $xs\not\in E(G)$, then we may also assume that $zs \not\in E(G)$ and $sy\in E(G)$. In this case, $H'$ contains the house graph induced by $\{x, y, z, a, b\}$ and yet $\deg_{H'}(y) \ge 4$, which cannot be. Therefore, we may assume that $H \in \mathcal{K}\cup \mathcal{W}\cup \mathcal{R}$. Assume first $H \in \mathcal{K}$. 
First, label the remaining vertex on the diamond as $w$ and assume that both $x$ and $y$ are the interior vertices of the diamond and that we have appended leaves $b_1, \dots, b_n$ to $w$. It follows that $\{s, w\} = S$ and $\{b_1, \dots, b_n, t\} = I$. Further, we may assume that $s$ is adjacent to at least one of $x$ or $z$. If $sx \in E(G)$, then both $\{sx, wy\}$ and $\{st, b_1w, yz\}$ are minimal edge dominating sets. Otherwise, if $sx \not\in E(G)$ and $sz \in E(G)$, then $\{sz, wy\}$ and $\{st, yz, wb_1\}$ are minimal edge dominating sets. Note that if $z$ is the exterior vertex of the diamond with appended leaves, then a similar argument as above works as well. Finally, assume $H \in \mathcal{W}\cup \mathcal{R}$. Note that if $H \in \mathcal{R}$, then $H$ is constructed from a windmill as $x$ and $y$ are assumed to have degree at least $3$ in $H$. Therefore, whether $H \in \mathcal{W}$ or $H \in \mathcal{R}$, we may assume that $H$ is constructed from a windmill $W$ which contains the house graph with vertices $\{x, y, z, a, b\}$ where $xa$ and $yb$ are edges, and $k\ge 0$ triangles  where the nose of the windmill is $z$. If $k\ge 1$, we label the remaining vertices of the windmill as $\{c_i, d_i: i \in [k]\}$ where $zc_id_i$ is a triangle in $H$ for $i \in [k]$. Moreover, if $H \in \mathcal{R}$, then said windmill has $n\ge 1$ $P_3$s attached to the windmill by identifying $z$ with one leaf from each of the $P_3$s. In this case, we label each of the remaining vertices of the $P_3$s as $\{e_i, f_i: i \in [n]\}$ where $e_if_iz$ is a path in $H$.  Note that if $s$ is adjacent to some vertex $u \in \{x, y\} \cup \{c_i, d_i: i \in [k]\}$, then $G - N_e[ab]$ contains the propeller with nose $z$ and yet $u$ does not have degree $2$, making it a smaller counterexample. Hence, we shall assume that the only neighbor of $s$ which is on a triangle is $z$. Consider $\widetilde{H} = G - N_e[ab]$ which we know is either $Cr$ or in $\mathcal{K}\cup \mathcal{P}\cup \mathcal{R}\cup \mathcal{W}$. $\widetilde{H}$ contains the propeller $W - \{a, b\}$ meaning that $\widetilde{H} \not\cong Cr$ and $\widetilde{H} \not\in \mathcal{K}\cup \mathcal{P}$. Moreover, $\widetilde{H} \not\in \mathcal{W}$ as $t$ is not adjacent to $z$. Thus, $\widetilde{H} \in \mathcal{R}$, $t$ is a leaf in $\widetilde{H}$, $N_G(s) \subseteq \{z, a, b, t\}$, and $N_G(t) \subseteq \{z, a, b, s\}$. If $W - \{a, b\}$ contains the triangle $c_1d_1z$, then $G - N_e[c_1d_1] \in \mathcal{R}$ and we may conclude that $G \in \mathcal{R}$. Similarly, if $\widetilde{H}$ contains the pendent edge $e_1f_1$, then $G- N_e[e_1f_1] \in \mathcal{R}$ and we may conclude $G \in \mathcal{R}$. However, now we have  $V(G) = \{x, y, z, a, b, s, t\}$ contradicting that our chosen counterexample has order at least $9$. 

Next, assume that exactly one vertex of $\{x, y, z\}$ has degree at least three in $H$. Without loss of generality, we may assume $\deg_H(x) \ge 3$. This implies that $H$ does not contain the house graph as an induced subgraph nor does it contain a diamond. Therefore,  $H \in \mathcal{P}\cup \mathcal{R}$ where $x$ is the nose of the propeller with $\deg_H(x) \ge 4$. Label the vertices of the remaining triangles other than $xyz$ (if they exist) as $xc_id_i$ for $i \in [k]$ and label the remaining $P_3s$ (if they exist) as $xe_jf_j$ for $j \in [m]$. Recall that $G$ contains at least two triangles. Therefore, if $H$ contains only one triangle, then we can interchange the roles of $st$ and $e_1f_1$ so that $H$ contains at least two triangles and $k \ge 1$. If $e_1f_1$ exists, then consider  $H'' = G-N_e[e_1f_1]$ where $\deg_{H''}(x) \ge 3$. If $\deg_{H''}(y) \ge 3$ or $\deg_{H''}(z)\ge 3$, then we have already considered this case above. So we may assume that $\deg_{H''}(y) = \deg_{H''}(z) = 2$ meaning that $N_G(s) \cap \{x, y, z\} = \{x\}$. On the other hand, if $e_1f_1$ does not exist, then $k \ge 2$ and we consider $G - N_e[c_1d_1]$  which again will lead us to $N_G(s) \cap \{x, y, z\} = \{x\}$. Lastly, $G - N_e[yz]$ contains the triangle $xc_1d_1$ and is therefore in $\mathcal{P}\cup \mathcal{R}$ where $x$ is the nose of the propeller as it is adjacent to $c_1, d_1$, and $s$. It follows that $G \in \mathcal{R}$.

\vskip2mm
\noindent\textbf{Case 2:} Suppose $I= \emptyset$.
\vskip2mm
We proceed by considering whether $S$ is an independent set or not. Suppose first that $S$ is independent. As we have assumed that $G$ contains at least two triangles, it follows that there is some vertex in $S$ that is adjacent to two vertices in $\{x, y, z\}$. Therefore, we assume that $s \in S$ is adjacent to both $x$ and $y$. If $\deg_G(z) = 2$, then every vertex in $S$ is adjacent to either $x$ or $y$ or both. In this case, $xy$ is an edge dominating set of $G$ and yet $\{xz, sy\}$ is a matching in $G$, which cannot be.  Therefore,  we may assume $\deg_G(z) \ge 3$. Let $s_z \in S - \{s\}$ be a neighbor of $z$ ($s$ is not adjacent to $z$ as $G$ is $K_4$-free). If there exists $s' \in S - \{s, s_z\}$ adjacent to $y$, then both $\{sx, yz\}$ and $\{sx, s'y, s_zz\}$ are maximal matchings in $G$. So we may assume $N_G(y) \subseteq \{x, z, s, s_z\}$. A similar argument shows that we may assume $N_G(x) \subseteq \{y, z, s, s_z\}$. If both $x$ and $y$ have degree $3$ in $G$, then $G \in \mathcal{K}$ and we are done. Thus, we shall assume either $xs_z\in E(G)$ or $ys_z\in E(G)$. Without loss of generality, we shall assume $xs_z \in E(G)$. As $G$ is $K_4$-free, it follows that $ys_z\not\in E(G)$. Further, $|S- \{s, s_z\}| \ge 4$ as $n(G) \ge 9$. Let $s' \in S- \{s, s_z\}$ which is adjacent to $z$. In this case, both $\{xs_z, ys, s'z\}$ and $\{xz, sy\}$ are maximal matchings in $G$. Hence, this case cannot occur and $S$ is not an independent set.

Therefore, we may assume that $S$ is not an independent set. Let $s_1, s_2 \in S$ be adjacent vertices in $G$. We consider $H'' = G- N_e[s_1s_2]$. Thus, $H'' \in \{Cr\} \cup \mathcal{K}\cup \mathcal{P}\cup \mathcal{R}\cup \mathcal{W}$. Note that we may assume $H'' \not\in \mathcal{R}$ as this would imply $I \ne \emptyset$. If $H''\cong Cr$, then we label the remaining vertices of the induced $P_4$ as $a, b, c$, and $d$ where $xa, yb$, and $zc$ are in $E(G)$. However, this implies that $d \in I$ which is a contradiction. Next, suppose $H'' \in \mathcal{K}$.  This implies that $H''$ is built from the diamond induced by $\{w, x, y, z\}$ where $w \in S$. Without loss of generality, we may assume $x$ and $y$ are the interior vertices of the diamond. Furthermore, since $w \in S$, it must be that $H''$ is built by attaching leaves $b_1, \dots, b_n$ to $z$ for otherwise $I \ne \emptyset$. If both $s_1$ and $s_2$ are adjacent to a vertex in $\{w, x, y\}$, then we could instead view $G$ as built from the triangle $wxy$ so that $z, s_1$, and $s_2$ are in $S$ and $I= \{b_1, \dots, b_n\}$, which is a case we have already considered. Therefore, we shall assume that $N_G(s_1) \cap \{x, y, z\} = \{z\}$ and $N_G(s_2) \cap \{x, y, z\} = \{z\}$. In this case, we could instead view $G$ as built from the triangle $zs_1s_2$ where $I = \{w\}$, yet another case that we have already considered. Therefore, $H'' \not\in \mathcal{K}$. 

Lastly, suppose $H'' \in \mathcal{P}\cup \mathcal{W}$ where $z$ is the nose of the propeller or windmill. Since $n(H'') \ge 7$, we may assume that $H''$ is either built from the house graph where $xyz$ is the triangle of the house graph and $k$ triangles of the form $zr_it_i$ for $i \in [k]$, or $H''$ is built from triangles $xyz$ and $zr_it_i$ for $i \in [k]$. In either case, $k \ge 1$. If $N_G(s_i) = \{z, s_j\}$ where $\{i,j\} = \{1, 2\}$, then $G \in \mathcal{P}\cup\mathcal{W}$.  Thus, we may assume, up to relabeling, that $s_1$ is adjacent to either $x, r_1$, or perhaps a vertex of degree two on the house graph. If $k \ge 2$, then $G - N_e[r_2t_2]$ contains the propeller or windmill induced by $V(G) - \{s_1, s_2, r_2, t_2\}$, meaning that $G - N_e[r_2t_2]\in \mathcal{P}\cup \mathcal{W}$ and $N_G(s_i) = \{s_j, z\}$ for $\{i, j\} = \{1, 2\}$, which is a contradiction. Thus, we shall assume that $V(H'') = \{x, y, z, a, b, r_1, t_1\}$ where $a$ and $b$ are the vertices of degree two on the house graph. If $s_1$ is adjacent to a vertex on the house graph other than $z$, then $G- N_e[r_1t_1]$ contains the house graph induced by $\{x, y, z, a, b\}$ and yet either $x$ or $a$ is adjacent to $s_1$, which does not fit the description  of any graph in $\{Cr, \mathcal{H}\}\cup \mathcal{K} \cup \mathcal{P} \cup \mathcal{R}\cup \mathcal{W}$. Thus, we may assume that $s_1$ (and similarly, $s_2$) is not adjacent to any vertex on the house graph. In this case, $G - N_e[ab]$ contains an induced propeller and it must be that $N_G(s_i) = \{s_j, z\}$ for $\{i, j\} = \{1, 2\}$, and $G \in \mathcal{W}$. Having exhausted all possibilities, no such counterexample exists. 

\end{proof}

We are now ready to prove we have completely characterized all well-edge-dominated graphs with girth $3$. Recall that $F_5$ is the fan of order $5$.

\begin{theorem} Let $G$ be $K_4$-free. $G$ is a connected, well-edge-dominated graph with girth $3$ if and only if $G \in \{ \mathcal{DH}, F_5,  Cr, W_1, W_2, W_3\}$ or $G \in \mathcal{G}$.
\end{theorem}

\begin{proof} Based on  Theorem~\ref{cor:G}, we only need to show that if $G$ is well-edge-dominated with girth $3$, then $G \in \mathcal{G} \cup \{\mathcal{DH}, F_5, Cr, W_1, W_2, W_3\}$. Moreover, if $G$ contains only one triangle, then we know by Theorem~\ref{thm:onetriangle} that $G \in \mathcal{T}\cup \mathcal{F} \cup \{K_3, Cr, \mathcal{H}, \mathcal{DH}\}$ and $\{\mathcal{H} ,K_3\} \cup \mathcal{T}\cup \mathcal{F} \subset \mathcal{G}$.  So we only consider when $G$ contains at least two triangles. We first note that if $n(G) \le 8$, then Sage verified that $G \in \mathcal{G} \cup \{W_1, W_2, W_3\}$. Therefore, we may assume that $n(G) \ge 9$. Assuming a counterexample exists, we choose the counterexample $G$ with the smallest order. We may choose a triangle $xyz$ and vertex $t$ in $G$ such that if we contract $x, y,$ and $z$ into a single vertex $T$, then   $\dist(T, t)$ is greatest among all choices of triangles $xyz$ and vertices $t$. Then among all neighbors of $t$ in $G$, choose $s$ so that $\dist(T,s)$ is maximum. Consider $G' = G - N_e[st]$. We may write $G' = G_0 \cup I$ where $I$ is an independent set in $G$ and $G_0$ is connected for otherwise $st$ was chosen incorrectly. We can partition $I = I_s \cup I_{st} \cup I_t$ where each vertex in $I_s$ is adjacent to $s$ and not $t$, each vertex in $I_t$ is adjacent to $t$ and not $s$, and each vertex in $I_{st}$ is adjacent to both $s$ and $t$. We claim that we may assume $I_t \cup I_{st} = \emptyset$. Suppose to the contrary that $\ell \in I_t \cup I_{st}$. If $\dist(T, s) < \dist(T, t)$, then $\dist(T, \ell) \ge \dist(T, t)$ which would imply that we should have chosen $s = \ell$. On the other hand, if $\dist(T, s) = \dist(T, t)$, then $\dist(T, \ell) > \dist(T, t)$ and we should have chosen $s=t$. In either case, we reach a contradiction. So we may assume that $I_t \cup I_{st} = \emptyset$. Furthermore, we may assume that $|I_s| \le 1$ for otherwise for some $\ell \in I_s$, $G - \ell$ is also well-edge-dominated, contradicting our choice of minimal counterexample. Thus, if $\dist(T, s) < \dist(T, t)$, then it is possible $\ell$ exists and is a leaf in $G$. Otherwise, if $\dist(T, s) = \dist(T, t)$, we should have chosen $\ell = t$ so this case cannot occur. It follows that $G' = G_0 \cup I_s$ where $|I_s|\le 1$ and $n(G_0)\ge 6$. We proceed by considering which graph $G_0$ is isomorphic to. Note that since $G$ is well-edge-dominated, $G_0$ is well-edge-dominated. Since $G$ is the smallest counterexample, $G_0 \in  \mathcal{G} \cup \{ \mathcal{DH}, F_5, Cr, W_1, W_2, W_3\}$. Since $n(G) \geq 9$, $G_0 \not\cong F_5$. We consider when $G_0  \in \{W_1, W_2, W_3\}$ (Case 1),  $G_0  \in \{\mathcal{DH}, Cr\}$ (Case 2), and $G_0 \in \mathcal{G}$ (Case 3).

\vskip2mm
\noindent\textbf{Case 1:} Suppose first that $G_0 \in \{W_1, W_2, W_3\}$. 
\vskip2mm
Label the vertices as in Figure~\ref{fig:Ws}. Suppose first that $\ell$ exists and $G$ has order $10$. Note that $H = G-N_e[ab]$ is well-edge-dominated with order $8$ and contains both a diamond and a leaf, namely $\ell$. Either $H$ is not connected or $H \in \{V_1, V_2, V_3\}$ in Figure~\ref{ND8}. Therefore, $s$ (and $t$) is not adjacent to any vertex in $\{w, x, y\}$. Similarly, $G - N_e[bc]$ is well-edge-dominated with order $8$ and therefore $s$ (and $t$) is not adjacent to any vertex in $\{x, y, z\}$. Thus, the only neighbors of $s$ in $G$ are in $\{a, b, c\}$. Consider $H' = G- N_e[cz]$ which contains the triangle $wxy$, no diamond, and leaf $\ell$. Suppose first that $H'$ is connected so that $H' \in \{V_4, V_5, \dots, V_{12}\}$ in Figure~\ref{ND8}. However, $H'$ contains at most two leaves and if $abs$ is a triangle, then $wxy$ and $abs$ (or $abt$) are vertex disjoint triangles with $wa\in E(G)$. Therefore, $H' \not\in \{V_4, V_5, V_{10}, V_{12}\}$. Furthermore, $\deg_{H'}(x) = 2 = \deg_{H'}(y)$ and $\deg_{H'}(w) = 3$ so $H' \not\in \{V_6, V_7, V_{11}\}$. Also, $\ell$ is not adjacent to $a$ so $H' \not\cong V_8$. It follows that $H' \cong V_9$ where $at \in E(G)$ and $bs \in E(G)$. However, this implies $G- N_e[at]$ is connected, with order $8$, contains a diamond and one leaf ($\ell$). The only graph in Figure~\ref{ND8} that satisfies these requirements is $V_3$, but $\ell$ is not adjacent to $z$. Therefore, this case cannot occur and we may assume that $H'$ is not connected, meaning that $s$ (and $t$) can only be adjacent to $c$ in $G_0$. Thus, $H$ is connected, with order $8$, contains a diamond and at most two leaves. Therefore, $H \in \{V_2, V_3\}$. In either case, this would imply $s$ is adjacent to $z$, which we have already ruled out as a possibility. Hence, we may assume that $\ell$ does not exist and $G$ has order $9$.

As before, setting $H = G- N_e[ab]$, we see that $H$ contains a diamond, but now has order $7$. Suppose first that $H$ is not connected. It follows that $s$ and $t$ are not adjacent to any vertex in $\{w, x, z, c\}$. Consider $H' = G-N_e[cz]$, which is a connected, well-edge-dominated graph of order $7$. By Figure~\ref{ND7}, no connected, well-edge-dominated graph of order $7$ contains three triangles except $W_{13}$. Thus, $\{a, b, s, t\}$ is not a diamond.  If $H'$ contains exactly one triangle, then it must be the case that $H' \cong W_9$ or $H'\cong W_{10}$ . In either case, $H'' = G - N_e[bc]$ is a connected, well-edge-dominated graph of order $7$ that contains a diamond; however, $H''$ is not $W_1$, $W_2$, $W_3$, or $W_4$, which is a contradiction. If $H''$ contains exactly two triangles, then those triangles do not share a common vertex, which cannot occur by Figure~\ref{ND7}. Thus, we may assume that $H$ is connected.  By Observation \ref{obs:diamondsize7}, $H \in \{W_1, W_2, W_3, W_4\}$.  It is not possible for $H$ to have three leaves. Thus, $H \not\cong W_4$ and we may assume that  $H \in \{W_1, W_2, W_3\}$. Hence,  $w$ is adjacent to either $s$ or $t$. In either case, $G - N_e[bc]$ is also in $\{W_1, W_2, W_3\}$ using the same reasoning that $H \in \{W_1, W_2, W_3\}$, and yet $w$ has degree four in $G-N_e[bc]$. Therefore, this case cannot occur either and we may assume $G_0 \not\in \{W_1, W_2, W_3\}$.

\smallqed

\vskip2mm
\noindent\textbf{Case 2:} Suppose that $G_0 \in \{\mathcal{DH}, Cr\}$.
\vskip2mm
Suppose first that $G_0 \cong\mathcal{DH}$ and label the vertices as in Figure~\ref{fig:houses} (a). Since $\alpha'(\mathcal{DH}) = 3$, it follows that $\alpha'(G) = 4$. Note that by our choice of $s$ and $t$, $s$ cannot be adjacent to $x, y$, or $z$ for otherwise we should have chosen $ab$ to act as the edge $st$. Moreover, by our choice of $s$ and $t$, this implies that we may assume $t$ is not adjacent to any vertex in $\{x, y, z\}$. Therefore, we have only to consider when $s$ is adjacent a vertex in $\{a, b, c, d\}$. Moreover, since $G$ is symmetric, we need only consider the two cases where $s$ is adjacent to $a$ versus the case that $s$ is adjacent to $c$.  

We consider $H'=G-N_e[bd]$ which is connected, has order $7$ or $8$ and contains the triangle $xyz$ where $\deg_{H'}(z) = \deg_{H'}(y)=2$ so $x$ is the nose of the triangle. If $H'$ contains only one triangle, then by Theorem~\ref{thm:onetriangle}, $\{s, t, a, c, x\}$ (or $\{s,t,a,c,x,\ell\}$ if $\ell$ exists) induces a bipartite graph $J$ in $G$ with bipartition $A_J \cup B_J$ where $|A_J|<|B_J|$ and $x$ is detachable. Suppose first that $\ell$ exists. This implies that $s \in A_J$, meaning $sa \in E(H')$ and $sc \not\in E(H')$. However, this implies that $x$ is not a detachable vertex in $J$. So we shall assume that $\ell$ does not exist. Moreover, we must have $s \in B_J$ meaning that $sc\in E(H')$ and $sa\not\in E(H')$. If $ta\not\in E(H')$, then again $x$ is not detachable in $J$. Now when we consider $G- N_e[xy]$, we have a well-edge-dominated (and thus equimatchable) graph induced by $\{s, t, a, b, c, d\}$ containing a perfect matching. Thus, by Theorem~\ref{thm:rmatch}, the graph induced by $\{s, t, a, b, c, d\}$ must be $K_{3, 3}$ so $sb$ and $td$ are edges in $G$. Now $\{sc, st, sb, dy, xz\}$ is a minimal edge dominating set contradicting the fact that $\alpha'(G)=4$. So we shall assume that $H'$ contains at least two triangles, and perhaps $\ell$ exists. Note that $H'$ is diamond-free for otherwise said diamond is induced by $\{s, t, a, c\}$ and the only graph in Figures~\ref{ND7} and \ref{ND8} that have an edge between a diamond and a triangle is $V_3$. However, $s$ is not adjacent to $x$ and $\ell$ is not adjacent to $c$ so this case cannot occur.  Moreover, the only diamond-free graphs in Figures~\ref{ND7} and \ref{ND8} that contain at least two triangles are in $\{W_8, W_{11}, W_{13}, V_{10}, V_{12}\}$. We know $H'$ cannot contain an induced propeller so the only possibility is $H' \cong V_{12}$, yet this is impossible as it implies $\ell$ is adjacent to $c$, which cannot be. Thus, we may assume $G_0 \not\cong\mathcal{DH}$.

Next, assume that $G_0 \cong Cr$ and label the vertices as in Figure~\ref{fig:houses} (c). If $s$ is adjacent to a vertex in $\{x, y, z\}$ and $\{t, d \}$ (or $\{t, d, \ell\}$ if $\ell$ exists) is an independent set, then $G$ satisfies the statement of Lemma~\ref{lem:small}. Thus, $G \in \{\mathcal{K} \cup \mathcal{P} \cup \mathcal{R}\cup \mathcal{W}\}  \subset \mathcal{G}$ since $\mathcal{K}\cup \mathcal{R} \subset \mathcal{G}$, and we are done. So we may assume that either $s$ is not adjacent to any vertex in $\{x, y, z\}$ or $td \in E(G)$. Suppose first that $s$ (and therefore $t$) is not adjacent to any vertex in $\{x, y, z\}$. As $G$ contains at least two triangles, we may assume one is $xyz$ and another contains $s$ and two vertices from $\{a, b, c, d\}$ or contains $s$ and $t$ and one vertex from $\{a, b, c, d\}$. Suppose first that either $abs, sta$, or $stb$ is a triangle in $G$. Then $G - N_e[cd]$ is connected, has order $7$ or $8$, contains two vertex disjoint triangles and yet $xyz$ is not on a diamond and each of $x$ and $y$ have degree $3$. No such graph exists in Figures~\ref{ND7} and \ref{ND8} so this case cannot occur. Similarly, if $bcs$ or $stc$ is a triangle in $G$, then $G- N_e[ad]$ satisfies the above and yet is not depicted in Figures~\ref{ND7} and \ref{ND8}. So a second triangle of $G$ must be one of $std, sad$, or $scd$. In this case, $G- N_e[by]$ contains a triangle and an induced $C_5$ ($\{a, d, c, x, z\}$) which is not depicted in Figures~\ref{ND7} and \ref{ND8}. Having exhausted all possibilities, we may assume that $s$ is adjacent to some vertex in $\{x, y, z\}$ and $td \in E(G)$. In this case, $G - N_e[td]$  is connected, has order $7$ or $8$, contains the triangle $xyz$ where each of $x, y$, and $z$ have degree at least three and at least one of $x, y$, or $z$ has degree $4$, and contains at most one leaf. Based on the graphs in  Figures~\ref{ND7} and \ref{ND8}, it follows that $G-N_e[cd] \in \{W_2, W_3, W_8, V_3\}$. However, $\{a, b, c\}$ induces a path so $G - N_e[cd] \not\in \{W_8, V_3\}$. It follows that $G- N_e[cd] \in \{W_2, W_3\}$, meaning that $\{s, x, y, z\}$ induces a diamond. In this case, exactly two vertices in $\{x, y, z\}$ have degree 4 in $G - N_e[cd]$, meaning $G - N_e[cd] \not\in \{W_2, W_3\}$, which is a contradiction. Therefore, $G_0 \not\cong Cr$. 

\smallqed

\vskip2mm
\noindent\textbf{Case 3:} Suppose that $G_0 \in \mathcal{G}$. 
\vskip2mm
We may assume that of $\{x, y, z\}$, $x$ has the maximum distance to $s$. Suppose first that $s$ (and therefore $t$) is not adjacent to any vertex of $\{x, y, z\}$. We proceed by considering whether $xyz$ is on a propeller, a windmill, or a diamond in $G_0$.
\begin{enumerate}
\item[1.] Suppose first that $xyz$ is on a propeller in $G_0$. We may assume that $z$ is the nose of the propeller. It follows that $\deg_{G_0}(x) = \deg_{G_0}(y)=2$. Consider $H = G-N_e[xy]$ which has order at least $7$. Since neither $s$ nor $t$ is adjacent to a vertex in $\{x, y, z\}$, we infer that $\deg_G(x) = \deg_G(y) = 2$. If $H \in \mathcal{G}$ where $z$ is the nose of a propeller, then $G \in \mathcal{G}$ and we are done.  Thus, we shall assume that either $H \in \{Cr, \mathcal{DH}, W_1, W_2, W_3\}$ or $H \in \mathcal{G}$ where $z$ is not the nose of a propeller. 

Suppose first that $H \in \{W_1, W_2, W_3\}$. This implies that $n(G)=9$ and $\ell$ does not exist. Label the vertices of the diamond as $\{a, b, c, d\}$ where $a$ and $c$ are the interior vertices. Further, if $H\cong W_2$, then assume $c$ is the interior vertex with degree $4$ and if $H\cong W_3$, then $a$ is the interior vertex with degree $4$ and b is adjacent to the vertex with degree $3$ that is not on the diamond. Label the remaining vertices of $H$ as $\{j, k, m\}$ where $j$ is adjacent to $d$ and $m$ is adjacent to $b$.  Note that $s, t$ and $z$ are all vertices in $H$. We proceed by considering each of the possibilities for $z$. Note that $\alpha'(G) = 4$. If $z = j$, then $\{jx, ac, bm\}$ is an edge dominating set of $G$, which is a contradiction. If $z= k$, then $\{kx, ad, bc\}$ is an edge dominating set of $G$. If $z = m$, then $\{xm, ac, dj\}$ is an edge dominating set of $G$. In each case, we reach a contradiction. So we may assume that $z$ is on the diamond in $H$. If $z = d$, then $\{xd, km, ac\}$ is an edge dominating set of $G$. If $z = a$, then $\{xa, cd, km\}$ is an edge dominating set of $G$. If $z = c$, then $\{cx, ab, j k\}$ is an edge dominating set of $G$. If $z = b$, then $\{bx, km, cd\}$ is an edge dominating set of $G$. In each case, we reach a contradiction so we may assume that $H \not\in \{W_1, W_2, W_3\}$.

Next, suppose that $H \cong \mathcal{DH}$. Again we may assume that $\ell$ does not exist. Label the vertices of the triangle in $H$ as $a, b$, and $c$ where $a$ has degree $2$ in $H$. Label the remaining vertices in $H$ with $j, k, m,$ and $r$ where $jkmr$ induces a path,  $b$ is adjacent to $j$, and $r$ is adjacent to $c$. Note that $s, t$, and $z$ are vertices in $H$. We consider the different possibilities for $z$ in $H$. Since $\mathcal{DH}$ is symmetric, we need only to consider when $z \in \{a, b, j, k\}$. Note that $\alpha'(G) = 4$. If $z = a$, then $st = km$ by choice of $s$ and $t$. This implies that  $\{jk, cr, xa\}$ is an edge dominating set  of $G$.  If $z = b$, then $\{bx, cr, km\}$ is an edge dominating set of $G$. If $z = j$, then $\{jx, mr, bc\}$ is an edge dominating set of $G$. If $z = k$, then $\{kx, cr, ab\}$ is an edge dominating set of $G$. In each case, we reach a contradiction so we may assume $H \not\cong \mathcal{DH}$. 

Next, assume $H \cong Cr$. As above, we may assume that $\ell$ does not exist. Label the vertices of the triangle in $H$ as $a, b, c$ and the remaining vertices in $H$ as $j, k, m$, and $r$ where $jkmr$ induces a path in $H$, $a$ is adjacent to $j$, $b$ is adjacent to $k$ and $c$ is adjacent to $m$. As above, we consider the different possibilities for $z$ in $H$. Further, as $Cr$ is symmetric, we only need to consider when $z \in \{a, b, j, k, r\}$. Note that $\alpha'(G) = 4$. If $z = a$, then $\{ax, bk, mr\}$ is an edge dominating set of $G$. If $z = b$, then $\{bx, aj, cm\}$ is an edge dominating set of $G$. If $z=j$, then $\{xj, km, ac\}$ is an edge dominating set of $G$. If $z=k$, then $\{xk, ac, mr\}$ is an edge dominating set of $G$. If $z = r$, then $\{rx, ac, bk\}$ is an edge dominating set of $G$. In each case, we reach a contradiction so this case cannot occur either. 

Finally, assume that $H\in \mathcal{G}$ yet $z$ is not the nose of a propeller in $H$. Suppose first that $H \in \mathcal{P}$ and built from $k$ triangles $T_1, \dots, T_k$ with $V(T_i) = \{a_i, b_i, c_i\}$ for $i \in [k]$. We may assume that $H$ is obtained from the disjoint union $T_1\cup \cdots \cup T_k$ by identifying the set $\{a_i: i \in [k]\}$ into a single vertex, call it $a$. Since we have assumed $z$ is not the nose of $P$, we may assume with no loss of generality that $z = b_1$. Moreover, $H$ contains $s$ and $t$ and $G_0 \in \mathcal{G}$. It follows that $a \not\in \{s, t\}$ and $st = b_ic_i$ for some $2 \le i \le k$. Further, since $G_0 \in \mathcal{G}$, $k=2$ as $z$ is the nose of the propeller containing $xyz$ in $G_0$. However, this implies $n(G) = 7$ which is a contradiction. Hence, we may assume $H \not\in \mathcal{P}$. 

Since $xyz$ is on a propeller in $G_0$, it is also the case that $H \notin \mathcal{W}$ and $H \not\cong K^*$. Moreover, $H$ was built from a connected, nontrivial, bipartite, and well-edge-dominated graph $J = (A_J\cup B_J, E_J)$ where $|A_J|<|B_J|$. Therefore, if $z$ is on a propeller in $H$, then $s$ and $t$ were chosen so that they are not vertices on the propeller with $z$ but instead vertices in $J$. However, this would imply that $G_0 \not\in \mathcal{G}$ as $z$ is the nose of its propeller in $G_0$ yet not the nose of the propeller in $H$, which is a contradiction. It follows that the induced propeller containing the triangle $xyz$ in $G_0$ is $K_3$ and $H$ contains another triangle as $G$ was assumed to contain at least two triangles. Now since $H$ is built from $J$, we may assume that we have chosen $B_{J}'\subset B_J$ and $Y\subseteq A_J$ such that 
 \begin{itemize}
\item $B_J'\cup Y \ne \emptyset$.
\item  $J- B_J'$ is a well-edge-dominated graph with $\gamma'(J-B_J') = \gamma'(J)$, $J-B_J'$ has no trivial components, and $|A_J| \le |B_J-B_J'|$. 
\item If $Y\ne \emptyset$, then we write $Y = \{y_1, \dots, y_{\ell}\}$ where $y_i$ is a support vertex in $J-B_J'$, and there are $\ell$ disjoint diamonds $D_1, \dots, D_{\ell}$ such that we identified an exterior vertex of $D_i$ with $y_i$. 
\end{itemize}
Moreover, if there exists an induced windmill in $H$, then we have chosen a nonempty $B_J'' \subseteq B_J'$ such that $B_J'' = \{s_1, \dots, s_k\}$ and there are $k$ disjoint windmills $W_1, \dots, W_k$ such that the nose of $W_i$ was identified with $s_i$. If $H$ contains an induced propeller (that doesn't contain the vertex we called $z$ from $xyz$), then we write $B_J' - B_J''= \{x_1, \dots, x_r\}$ and there are $r$ disjoint propellers $P_1, \dots, P_r$ such that the nose of $P_i$ was identified with $x_i$. Additionally, every neighbor of $s_i$ (should it exist)  in $A_J$ is a support vertex in $J-B_J'$. Now, because $J$ is nontrivial and $|A_J|<|B_J|$,  we know that $|A_J|+|B_J| \ge 3$.  If we can show that $\widetilde{B_J} = \{x_0, x_1, \dots, x_r, s_1, \dots, s_k\}$ is such that $J - \widetilde{B_J}$ is well-edge-dominated with no trivial components where $x_0$ is the vertex in $J$ that we identified $z$ with to create $G_0$ and every neighbor of $s_i$ in $A_J$ is a support vertex in $J-\widetilde{B_J}$, then $G \in \mathcal{G}$.  We consider the diameter of $J$. Suppose first that $J$ is a star. It follows that there is some vertex $v \in B_J - B_J'$ for otherwise $J - B_J'$ contains a trivial component. Therefore, $x_0 \in B_J-B_J'$ and if $|B_J - B_J'| \ge 2$, then $G \in \mathcal{G}$. So we shall assume that $B_J - B_J' = \{x_0\}$. Thus, $B_J = \{x_0, x_1, \dots, x_r, s_1, \dots, s_k\}$ and $A_J = \{a\}$. If $P_1$ exists with triangle $jkm$ such that $j$ is the nose of $P_1$, then the component $\widetilde{G}$ of $G-N_e[jk]$ containing $a$ is not well-edge-dominated as it is not in $\{Cr, \mathcal{DH}, W_1, W_2, W_3\}$ nor is it in $\mathcal{G}$. We also reach a similar contradiction if $W_1$ exists. However, this would imply that $J = K_2$ and $|B_J| = |A_J|$ which is another contradiction. 

Thus, we may assume that $J$ is not a star. This implies $|A_J| \ge 2$ and $|B_J| \ge 3$. Consider $H' = G-N_e[xz]$. Since $y$ is on a propeller of $G_0$ and it is assumed that $y$ is not adjacent to $s$ or $t$, it follows that $y$ is an isolate in $H'$. Moreover, $H'$ contains a triangle as $G$ contains at least two triangles. Therefore, we consider $H_1'$, a component of $H'$ of maximum order amongst all components containing a triangle. We first claim that we may assume that $H' = H_1' \cup \{y\}$. To see this note that we originally chose $s$ and $t$ to have maximum distance from the triangle $xyz$. Therefore, if $H' - \{y\}$ contains at least two components, say $H_1'$ and $H_2'$, then $H_1'$ contains a triangle and $st$ is in $H_2'$. However, in this case, any triangle in $H_1'$ is further from $st$ than $xyz$, which is a contradiction. Therefore, we may assume that $H_1' = H' - \{y\}$. Suppose first that  $H_1'$ contains an induced windmill. Hence, $n(H_1') \ge 8$ and we have $H_1' \in \mathcal{G}$. Thus, $H_1'$ is obtained from the connected, nontrivial, bipartite and well-edge-dominated graph $\widetilde{J} = (A_{\widetilde{J}} \cup B_{\widetilde{J}}, E_{\widetilde{J}})$ where  $\widetilde{J} = J - \{x_0\}$, $B_{\widetilde{J}} = B_J - \{x_0\}$, $A_{\widetilde{J}} = A_J$, and $|A_{\widetilde{J}}|<|B_{\widetilde{J}}|$. Moreover, $x_0 \not\in B_J'$ so we may assume that we have chosen $B_{\widetilde{J}}' = B_J'$ to create $H_1'$. It follows that  $\widetilde{J} - B_{\widetilde{J}}'$ is well-edge-dominated as $H_1' \in \mathcal{G}$. Furthermore, $J - \widetilde{B_J} \cong \widetilde{J} - B_{\widetilde{J}}'$, meaning $J - \widetilde{B_J}$ is well-edge-dominated, and each neighbor of $s_i$ in $A_J$ is a support vertex in $\widetilde{J} - B_{\widetilde{J}}'$.  Hence,  $G \in \mathcal{G}$. Similarly, if $H_1'$ contains two induced diamonds, or $H_1'$ contains an induced diamond and an induced propeller, then $n(H_1') \ge 8$ and both $H_1'$ and $G$ are in $\mathcal{G}$. On the other hand, since $H \in \mathcal{G}$,  every triangle in $H$ (and therefore $G$ and $H_1'$) is either on a diamond where an exterior vertex has degree two in $H$, or has exactly one vertex with degree at least three in $H$. Hence, since $n(H_1') \ge 6$ and $H'_1 \not\in \{W_1, W_2, W_3, \mathcal{DH}, Cr\}$, it must be that $H_1' \in \mathcal{G}$. Using the same argument as above, we have that $G \in \mathcal{G}$.

%%%%%%%%%%%%%%%%%%%%%%%%%%%%%%%%%%%%%%%%%%%%%%%%%%%%%%%%%%
\item[2.] Next, suppose $xyz$ is on an induced windmill $W$ in $G_0$ where $z$ is the nose of the windmill.  Suppose first that $G_0 \in \mathcal{W}$ built from $\mathcal{H}$ with vertices $\{a, b_1, c_1, d, e\}$ where $ab_1c_1$ is the triangle, $b_1d$ is an edge, and all remaining triangles in $G_0$ are of the form $ab_ic_i$ for $2 \le i \le k$ (note that $k\ge 2$ as $n(G_0) \ge 7$). Since neither $s$ (nor $t$) is adjacent to any vertex in $\{x, y, z\}$, we may assume that $s$ is not adjacent to $a$ in $G$. Suppose first that $s$ is adjacent to $b_i$ (or $c_i$) for some $i \in [k]$. Then $G- N_e[de]$ is not well-edge-dominated as there is no graph in $\mathcal{G} \cup \{Cr, \mathcal{DH}, W_1, W_2, W_3\}$ containing an induced propeller where a vertex other than the nose of the propeller has degree at least three in the graph. So the only neighbors of $s$ (and therefore $t$) in $G$ are in $\{d, e\}$. Without loss of generality, assume that $sd \in E(G)$.  Let $H' = G- N_e[b_2c_2]$ which is connected and $n(H') \ge 7$. As $H'$ contains the house graph as an induced subgraph (with vertices $\{a, b_1, c_1, d, e\}$) and $\deg_{H'}(c_1) = \deg_{H'}(b_1) = 3$ and $\deg_{H'}(d) \ge 3$, it follows that $H' \cong \mathcal{DH}$. In this case, both $\{st, de\} \cup \{b_ic_i: i \in [k]\}$ and $\{et, b_1d\} \cup \{ab_i: 2 \le i \le k\}$ are minimal edge dominating sets in $G$, another contradiction. Thus, we may assume $G_0 \not\in \mathcal{W}$. 

 %Let $\{x, y, z, a, b\}$ be the vertices on the house graph in $W$ where $xa$ and $yb$ are edges in $W$. 
Therefore, we can construct $G_0$ from the connected, nontrivial, bipartite, and well-edge-dominated graph $J = (A_J\cup B_J, E_J)$ where $|A_J|<|B_J|$ such that we have chosen $B_J'\subset B_J$ and $Y \subseteq A_J$ where  \begin{itemize}
\item $B_J'\cup Y \ne \emptyset$.
\item  $J- B_J'$ is a well-edge-dominated graph with $\gamma'(J-B_J') = \gamma'(J)$, $J-B_J'$ has no trivial components, and $|A_J| \le |B_J-B_J'|$. 
\item If $Y\ne \emptyset$, then we write $Y = \{y_1, \dots, y_{\ell}\}$ where $y_i$ is a support vertex in $J- B_J'$ and there are $\ell$ disjoint diamonds $D_1, \dots, D_{\ell}$ such that an exterior vertex of $D_i$ is identified with $y_i$.  
\end{itemize}
Note that we may assume $B_J''\subseteq B_J'$ is nonempty and we can write $B_J'' = \{s_1, \dots, s_k\}$ such that for $k$ disjoint windmills $W_1, \dots, W_k$ (where $W=W_1$), we have identified the nose of $W_i$ with $s_i$ and every neighbor of $s_i$ in $A_J$ is a support vertex in $J - B_J'$. Moreover, if $B_J' - B_J'' \ne \emptyset$, we enumerate $B_J' - B_J'' = \{x_1, \dots, x_r\}$ such that there are $r$ disjoint propellers $P_1, \dots, P_r$ where the nose of $P_i$ is identified with $x_i$. We know that $|A_J|+|B_J| \ge 3$ as there is some neighbor $u \in A_J$ adjacent to $s_1$ which is a support vertex in $J-B_J'$. We let $v$ be a leaf of $u$ in $J - B_J'$. Regardless of whether $xyz$ is on the house graph in $W_1$ or not, we label the vertices of the house graph in $W_1$ as $\{z, a, b, c, d\}$ where $zcd$ is the triangle and $da \in E(G)$. We argue that $\deg_G(a) = \deg_G(b) = 2$. Suppose to the contrary this is not the case. Since $\deg_{G_0}(a) = \deg_{G_0}(b) = 2$, it must be that $s$ is adjacent to either $a$ or $b$ in $G$. Without loss of generality, assume $sa \in E(G)$. Let $\widetilde{G}$ be the component of $G-N_e[uv]$ containing $z$. Therefore, $\widetilde{G}$ contains an induced windmill with nose $z$ yet the degree of $a$ is at least three. Since $z$ is not adjacent to $s$ or $t$, it follows that $\widetilde{G} \cong \mathcal{DH}$ and $W_1 \cong \mathcal{H}$. Thus, $zcd = zxy$, $\dist(T, s) = 2$ and since $\widetilde{G} \cong \mathcal{DH}$, neither $s$ not $t$ is adjacent to any vertex in $G$ other than potentially $\{u, v, a, b\}$. We know that $G$ contains a second triangle and if that triangle is on $P_1$ or $W_2$ in $G_0$, then $st$ was chosen incorrectly. Therefore, $B_J' = \{s_1\}$. Similarly, $st$ was chosen incorrectly unless $J$ is a star and we can write $A_J = \{u\}$ and $B_J = \{v_1, \dots, v_m, s_1\}$ where $v_1 = v$. If $D_1$ exists with vertices $\{u, w_1, w_2, w_3\}$ where $w_3$ is an exterior vertex, then we should have chosen $xyz = w_1w_2w_3$ as neither $s$ nor $t$ is adjacent to any vertex on this triangle. Hence, $D_1$ does not exist either and the second triangle in $G$ involves $st$. In this case, $G-N_e[ab]$ is connected and in $\mathcal{G}$ with triangle $stv_1$ and $v_2$ exists. But then $\{as, tb, ad, cb, uv_2\}$ is a minimal edge dominating set of $G$ as well as $\{st, ab, cd, uv_2\}$, the final contradiction. Thus, $\deg_G(a) = \deg_G(b) = 2$.

Now consider $H =G-N_e[ab]$ where $xyz$ is on a propeller with nose $z$ and $H\in \mathcal{G}$ as $H \not\in \{\mathcal{DH}, Cr, W_1, W_2, W_3\}$. Moreover, $H$ contains at least two triangles. Thus, $H$ is built from the connected, nontrivial, bipartite, and well-edge-dominated graph $\widetilde{J} = (A_{\widetilde{J}}\cup B_{\widetilde{J}}, E_{\widetilde{J}})$ where $|A_{\widetilde{J}}|<|B_{\widetilde{J}}|$ such that we have chosen $B_{\widetilde{J}}' \subset B_{\widetilde{J}}$ and $\widetilde{Y} \subseteq A_{\widetilde{J}}$ where  \begin{itemize}
\item $B_{\widetilde{J}}'\cup \widetilde{Y} \ne \emptyset$.
\item  $\widetilde{J} - B_{\widetilde{J}}'$ is a well-edge-dominated graph with $\gamma'(\widetilde{J} - B_{\widetilde{J}}') = \gamma'(\widetilde{J})$, $\widetilde{J} - B_{\widetilde{J}}'$ has no trivial components, and $|A_{\widetilde{J}}| \le |B_{\widetilde{J}} - B_{\widetilde{J}}'|$. 
\item If $\widetilde{Y}\ne \emptyset$, then we write $\widetilde{Y} = \{t_1, \dots, t_{\ell'}\}$ where $t_i$ is a support vertex in $\widetilde{J}-B_{\widetilde{J}}'$. 
\end{itemize}
Notice that $J = \widetilde{J} - \{s, t\}$ (if either is in $J$), and $H$ contains the propeller with triangle $xyz$. If $H$ contains an induced windmill, then some $B_{\widetilde{J}}'' \subseteq B_{\widetilde{J}}'$ has been chosen so that each vertex in $B_{\widetilde{J}}''$ has been identified with a nose of exactly one windmill. Moreover, $B_{\widetilde{J}}' - B_{\widetilde{J}}'' \ne \emptyset$ as the propeller containing $xyz$ has been identified with $s_1 \in B_{\widetilde{J}}' - B_{\widetilde{J}}''$ in $H$. Therefore, regardless of whether $st$ is on a diamond, a windmill, a propeller, or in $\widetilde{J}$, all we must show is that every neighbor of $s_1$ in $A_{\widetilde{J}}$ is a support vertex in $\widetilde{J} - B_{\widetilde{J}}'$.    Suppose to the contrary that there is some $u' \in N_{\widetilde{J}}(s_1)$ that is not a support vertex in $\widetilde{J} - B_{\widetilde{J}}'$. Note that $u' \not\in \{s, t\}$ as neither $s$ nor $t$ is adjacent to $z$. As $u'$ was a support vertex in $J - B_J'$ with leaf $v'$, this means that in $H$ $v'$ is adjacent to either $s$ or $t$. Without loss of generality, let us assume $sv'\in E(G)$. Note first that if the propeller in $H$ containing $zcd$ is of order at least $6$, then we may assume that $ztw$ is also on a triangle on the propeller as well as the windmill $W_1$ in $G_0$. Then $G - N_e[tw]$ is in $\mathcal{G}$ and built from the same graph $\widetilde{J}$ as $H$ and same choice of $B_{\widetilde{J}}'\cup \widetilde{Y}$. However, $B_{\widetilde{J}}''$ differs with respect to $G - N_e[tw]$ than with respect to $H$ as $z$ will be the nose of a windmill in $G- N_e[tw]$ rather than the nose of a propeller in $H$. However, every neighbor of $s_1$ in $A_{\widetilde{J}}$ is a support vertex in $\widetilde{J} - B_{\widetilde{J}}'$ as $G - N_e[tw]\in \mathcal{G}$, contradicting the assumption that $u'$ is not a support vertex in $H$. Therefore, we shall assume the propeller in $H$ containing $zcd$ is $K_3$. Since $H$ contains at least two triangles, there exists some triangle $klm$ on a windmill, propeller, or diamond in $H$, where we shall assume $k$ is the nose. If $klm$ is on a windmill or propeller, then $G - N_e[kl]$ contains the component $\overline{H}$ with windmill induced by $\{z, a, b, c, d\}$ as well as the vertices $\{u', v', s, t\}$. Thus, $\overline{H}$ has order at least $9$ and is in $\mathcal{G}$. Therefore, it is built from a connected, nontrivial, bipartite, and well-edge-dominated graph $T= (A_T\cup B_T, E_T)$ where $|A_T|<|B_T|$ with $s_1 \in B_T$. Furthermore, $N_{\widetilde{J}}(s_1) = N_T(s_1)$ and yet $u'$ will be a support vertex in $\overline{H}$ by assumption that $\overline{H} \in \mathcal{G}$. However, this cannot be as it would imply that $u'$ is a support vertex in $H$. Thus, we shall assume that the only other triangles in $H$ are on diamonds. Now if $k$ was identified with $u'$ to create $H$, then by definition of $\mathcal{G}$, $u'$ is a support vertex in $H$, which is another contradiction. On the other hand, $G- N_e[lm]$ has the component $\widetilde{H}$ containing the windmill induced by $\{z, a, b, c, d\}$ as well as vertices $\{u', v', s, t\}$ and is therefore in $\mathcal{G}$. However, $u'$ is a support vertex in $\widetilde{H}$, implying $u'$ is a support vertex in $H$, the final contradiction. Having exhausted all cases, we may conclude that no such $u'$ exists and $G \in \mathcal{G}$.

%Let $F$ be any minimal edge dominating set of $G_0$. We can also assume that the remaining vertices in $W_1$ other than those on the house graph induced by $\{z, a, b, c, d\}$ are of the form $b_i, c_i$ for $i\in [k]$ where $zb_ic_i$ is a triangle. Note that $F$ contains an edge incident to $u'$ so $|F\cap E(W_1)| =k+2$. Let $F' = F \cap E(W_1)$ together with the edge in $F$ incident to $u'$.  One can easily verify that 
%\[(F-F') \cup \{sv'\} \cup \{b_ic_i: i \in [k]\} \cup \{zc, ad\}\] is a minimal edge dominating set of $G$ with cardinality 
%\[|F-F'| +1 + k+2= |F|.\]
%On the other hand, $F\cup \{st\}$ is also a minimal edge dominating set of $G$, which is a contradiction. It follows that no such $u'$ exists and $G \in \mathcal{G}$. 

\item[3.] Suppose that $xyz$ is on a diamond and $w$ is the fourth vertex of the diamond in $G_0$. Without loss of generality, we may assume $x$ and $y$ are interior vertices of the diamond. We shall assume first that $w$ is the nose of the diamond and consider $H = G - N_e[yz]$ which is connected and has order at least $7$ as $s$ is not adjacent to either $y$ or $z$. Moreover, $x$ is a leaf of $w$ in $H$ so that $H \in \mathcal{G}$ and $H \not\in \mathcal{P}\cup \mathcal{W}$. This implies that $H$ is built from a connected, nontrivial, bipartite, and well-edge-dominated graph $J = (A_J \cup B_J, E_J)$ where $|A_J| < |B_J|$ and we have chosen $B_J' \subset B_J$ and $Y\subseteq A_J$ such that \begin{itemize}
\item $B_J'\cup Y \ne \emptyset$.
\item  $J- B_J'$ is a well-edge-dominated graph with $\gamma'(J-B_J') = \gamma'(J)$, $J-B_J'$ has no trivial components, and $|A_J| \le |B_J-B_J'|$. 
\item If $Y\ne \emptyset$, then we write $Y = \{y_1, \dots, y_{\ell}\}$ where $y_i$ is a support vertex in $J- B_J'$. 
\end{itemize}

 Furthermore, $x \in B_J$ and $w \in A_J$. Set $J' = J - \{x\}$, $A_{J'} = A_J$, $B_{J'} = B_J - \{x\}$, and $B_{J'}' = B_J'$. If we can show that $J' =(A_{J'} \cup B_{J'}, E_{J'})$ is a connected, well-edge-dominated bipartite graph with $|A_{J'}| < |B_{J'}|$ and $J' - B_{J'}' = J' - B_J'$ is well-edge-dominated with no trivial components such that $w$ is still a support vertex in $J' - B_J'$, then we have shown that $G \in \mathcal{G}$. First, consider $H' = G- N_e[xy]$ which is also connected, but has order at least $6$ as $z$ is an isolate in $H'$. However, since $H \not\in \{Cr, \mathcal{DH}, W_1, W_2, W_3\}$, $H' \not\in \{Cr, \mathcal{DH}, W_1, W_2, W_3\}$ so $H' \in \mathcal{G}$. Thus, since $H$ and $H'$ differ only by a leaf, this means that we can build $H'$ from  $J'$ where $J'$ (the same defined above) is connected, nontrivial, bipartite, and well-edge-dominated with $|A_{J'}| < |B_{J'}|$. Moreover, to build $H'$ it must be that we chose $B_{J'}' = B_J'$ and by assumption that $H' \in \mathcal{G}$,  $J' - B_J'$ is well-edge-dominated with no trivial components. We only need to show that $w$ is a support vertex in $J' - B_J'$. Suppose for the sake of contradiction that $w$ is not a support vertex in $J' - B_J'$. We construct two different minimal edge dominating sets for $G$. If $H$ was constructed using propellers, then enumerate those propellers as $P_1, \dots, P_r$ and for each $P_i$, let $E_i$ be the set of edges which are not incident to the nose of $P_i$. If $H$ was constructed using windmills, then enumerate those windmills as $W_1, \dots, W_k$ and for each $W_j$, let $F_j$ be any maximal matching of $W_j$ with edges that are not incident to the nose of $W_j$. Finally, if $H$ was constructed using diamonds, then enumerate those diamonds as $D_1, \dots, D_{\ell}$ and for each $D_i$, let $e_i$ be the edge between two interior vertices of $D_i$. Since $w$ is not a support vertex in $J'$, by Lemma~\ref{lem:bipartitesupports}, there exists a minimal edge dominating set $\widetilde{F}$ of $J'$ that does not contain an edge incident to $w$. Then both 
\[\widetilde{F} \cup \bigcup_{i=1}^r E_i \cup \bigcup_{j=1}^kF_j \cup \{e_1, \dots, e_{\ell}\} \cup \{xy\}\]
and 
\[\widetilde{F} \cup \bigcup_{i=1}^r E_i \cup \bigcup_{j=1}^kF_j \cup \{e_1, \dots, e_{\ell}\} \cup \{xz, yz\}\]
are minimal edge dominating sets of $G$, which cannot be. Thus, $w$ is a support vertex in $J'$ and $G \in \mathcal{G}$. 

Lastly, assume that $z$ is the nose of the diamond containing $xyz$ and $\deg_{G_0}(w) = 2$. If $w$ is not adjacent to $s$ or $t$, then we should have chosen $wxy$ to act as the triangle where $\dist(T, t)$ is maximum. Therefore, we shall assume that $sw \in E(G)$. It follows that $\dist(T, s) = 2$ and $\dist(T, t) \le 3$. Moreover, $G_0$ is built from a connected, nontrivial, bipartite, and well-edge-dominated graph $\widetilde{J}= (A_{\widetilde{J}}\cup B_{\widetilde{J}}, E_{\widetilde{J}})$ where $|A_{\widetilde{J}}|<|B_{\widetilde{J}}|$ and we have chosen $B_{\widetilde{J}}' \subset B_{\widetilde{J}}$ and $\widetilde{Y}\subseteq A_{\widetilde{J}}$ 
\begin{itemize}
\item $B_{\widetilde{J}}'\cup \widetilde{Y} \ne \emptyset$.
\item  $\widetilde{J} - B_{\widetilde{J}}'$ is a well-edge-dominated graph with $\gamma'(\widetilde{J} - B_{\widetilde{J}}') = \gamma'(\widetilde{J})$, $\widetilde{J} - B_{\widetilde{J}}'$ has no trivial components, and $|A_{\widetilde{J}}| \le |B_{\widetilde{J}}-B_{\widetilde{J}}'|$. 
\item If $\widetilde{Y}\ne \emptyset$, then we write $\widetilde{Y} = \{y_1, \dots, y_{\ell}\}$ where $y_i$ is a support vertex in $\widetilde{J} - B_{\widetilde{J}}'$. 
\end{itemize}

 Furthermore, we have identified $z$ with $a \in A_{\widetilde{J}}$ which is a support vertex in $\widetilde{J} - B_{\widetilde{J}}'$ with leaf $b_{z}$. Suppose there exists $a' \in A_{\widetilde{J}} - \{a\}$. As $\widetilde{J}$ is connected, we may assume that we chose $a'$ such that $a$ and $a'$ share a common neighbor, say $b'$, in $B_{\widetilde{J}}$. It follows that the component  $\widetilde{H}$ in $G - N_e[a'b']$ containing $w$ contains the vertices $\{w, x, y, z, s, t, b_z\}$. Thus, $\widetilde{H} \in \{Cr, \mathcal{DH}, W_1, W_2, W_3\}\cup \mathcal{G}$. However, $\widetilde{H}$ contains a diamond and therefore $\widetilde{H} \not\in \{Cr, \mathcal{DH}\}$. On the other hand, $\widetilde{H}$ contains a diamond where both exterior vertices have degree at least three and neither $s$ nor $t$ is adjacent to $x$ or $y$. It follows that $\widetilde{H} \cong W_1$ and $tb_z \in E(G)$. Now when we consider $G - N_e[tb_z]$, we have a connected graph with order at least seven and a diamond where both exterior vertices have degree at least three. Thus, $G - N_e[tb_z] \cong W_1$ meaning that $\widetilde{J} \cong P_4$, which is not a well-edge-dominated bipartite graph. It follows that no such $a' \in A_{\widetilde{J}} - \{a\}$ exists and $\widetilde{J}$ is a star and we can write $B_{\widetilde{J}} = \{b_1, \dots, b_n\}$ where $b_n = b_z$. If $t$ is adjacent to some $b_i$, then $G- N_e[tb_i]$ is connected with order at least seven and contains a diamond where both exterior vertices have degree at least three. Thus, $G - N_e[tb_i] \in \{W_1, W_2, W_3\}$ and yet this cannot be as $B_{\widetilde{J}}$ is an independent set. Therefore, this case cannot occur either and we may conclude that $G \in \mathcal{G}$. 

\end{enumerate}

Finally, we assume that $s$ is adjacent to at least one vertex of $\{x, y, z\}$. In this case, we can write $V(G) = \{x, y, z\} \cup S \cup T$ where each vertex in $S$ has a neighbor in $\{x, y, z\}$, and each vertex in $T$ has a neighbor in $S$, but no neighbor in $\{x, y, z\}$, and $T$ is an independent set in $G$. However, now $G$ satisfies the statement of Lemma~\ref{lem:small}, meaning that $G \in \{Cr, \mathcal{H}\}\cup \mathcal{K} \cup \mathcal{P} \cup \mathcal{R}\cup \mathcal{W}$. Since $\mathcal{K}\cup \mathcal{R} \subset \mathcal{G}$, we are done. 

\end{proof}

\section{Concluding Remarks} 
We now have a characterization for all non-bipartite, well-edge-dominated graphs that do not contain an induced $K_4$. Therefore, the final step to characterizing all non-bipartite, well-edge-dominated graphs would be to characterize those that contain an induced $K_4$. Note that there are instances when one can create a well-edge-dominated graph containing a $K_4$ from a graph in $\mathcal{G}$. For instance, if we take $K^* \in \mathcal{G}$ and add the missing edge on the diamond, we end up with a well-edge-dominated graph containing a $K_4$. So we pose the following question:
\vskip5mm

\noindent\textbf{Problem 1:} Is it possible that all non-bipartite well-edge-dominated graphs containing a $K_4$ are the result of adding missing edges from diamonds of a graph in $\mathcal{G}$?

\vskip5mm

Next, we point out that ideally one would be able to characterize all equimatchable graphs. Is it possible that a characterization for $K_4$-free equimatchable graphs is attainable? We pose this as another open problem:

\vskip5mm

\noindent\textbf{Problem 2:} Characterize all connected, non-bipartite, $K_4$-free, equimatchable graphs.

\section{Conflict of Interest}

The authors declare that they have no conflicts of interest.

\section{Appendix}

\subsection{Classes of Graphs}

\begin{itemize}
\item \textbf{Class $\mathcal{F}$}: A graph $G \in \mathcal{F}$ if $G$ is obtained from the disjoint union of the house graph $\mathcal{H}$, depicted in Figure~\ref{fig:houses}(b), and a well-edge-dominated bipartite graph $G'=( A \cup B, E)$ where  $|A|< |B|$,  by identifying the vertex of degree two in $\mathcal{H}$ on the triangle with $w \in V(G')$ where $w$ is strongly detachable. 
\item \textbf{Class $\mathcal{T}$}: A graph $G \in \mathcal{T}$ if $G$ is obtained from the disjoint union of $K_3$ and a well-edge-dominated bipartite graph $G'=(A \cup B, E)$ where  $|A| < |B|$,  by identifying $x \in V(K_3)$ with $w \in V(G')$ where $w$ is detachable.
\item \textbf{Class $\mathcal{P}$ (Propellers)}: Let $T_1, \dots, T_k$ be $k$ disjoint triangles and for each $i \in [k]$ label one vertex of triangle $T_i$ by $v_i$. A graph $G \in \mathcal{P}$ if it is obtained from $T_1, \dots, T_k$ by identifying each vertex in $\{v_1, \dots, v_k\}$ into a single vertex.
\item \textbf{Class $\mathcal{W}$ (Windmills)}: Let $T_1, \dots, T_k$ be $k$ disjoint triangles and for each $i \in [k]$ label one vertex of triangle $T_i$ by $v_i$. Let $\mathcal{H}$ be the house graph, depicted in Figure~\ref{fig:houses}(b), and label the vertex on the triangle with degree two $x$. A graph $G \in \mathcal{W}$ if it is obtained from $T_1, \dots, T_k,$ and $ \mathcal{H}$ by identifying each vertex in $\{v_1, \dots, v_k, x\}$ into a single vertex.  
\item \textbf{Class $\mathcal{G}$}: A graph $G \in \mathcal{G}$ if:
\begin{itemize}
\item The graph $G\cong K^*$, where $K^*$ is obtained from the disjoint union of a diamond and $K_2$ by identifying a vertex of $K_2$ with an exterior vertex of the diamond.
\item $G \in \mathcal{P}\cup \mathcal{W}$.
\item $G$ can be constructed in the following way. Let $G'=(A \cup B, E')$ be a connected, nontrivial, bipartite and well-edge-dominated graph with $|A| < |B|$.  Choose a set $B' \subset B$ and a set $Y \subseteq A$ such that the following are true:
\begin{itemize}
\item $B'\cup Y \ne \emptyset$.
\item  $G'- B'$ is a well-edge-dominated graph with $\gamma'(G'-B') = \gamma'(G')$, $G'-B'$ has no trivial components, and $|A| \le |B-B'|$. 
\item If $Y\ne \emptyset$, then we write $Y = \{y_1, \dots, y_{\ell}\}$ where $y_i$ is a support vertex in $G'- B'$. 
\end{itemize}
 Each vertex $v$ in $B'$ is called detachable and it is called strongly detachable if each neighbor of it in $G'$ is a support vertex in $G'-B'$.  If $Y\ne \emptyset$, let $D_1, \dots, D_{\ell}$ be a disjoint union of diamonds, and for each $D_i$, label one of the vertices of degree two as the ``nose" of the diamond. If $B' \ne\emptyset$, then choose a set $B'' \subseteq B'$ (possibly empty) where each vertex in $B''$ is strongly detachable. If $B'' \ne \emptyset$, then write $B'' = \{s_1, \dots, s_k\}$ and choose any set of $k$ windmills, enumerated as $W_1, \dots, W_k$. If $B' - B'' \ne \emptyset$, write $B' - B'' = \{x_1, \dots, x_r\}$ and choose any set of $r$ propellers, enumerated as $P_1, \dots, P_r$.  We obtain $G$ from $G'$ under the following rules: 
\begin{enumerate}
\item[(a)] The nose of $W_i$ is identified with $s_i$ in $B''$.
\item[(b)] The nose of $P_i$ is identified with $x_i$ in $B' - B''$. 
\item[(c)] The nose of $D_i$ is identified with the support vertex $y_i$ of $G'$ (which stays a support vertex in $G$). 
\end{enumerate}
\end{itemize}

\item \textbf{Class $\mathcal{K}$ (Kites)}: A graph $G \in \mathcal{K}$ if $G$ is obtained by appending $n\ge 1$ leaves to exactly one exterior vertex of a diamond. 
\item \textbf{Class $\mathcal{R}$}: A graph $G \in \mathcal{R}$ if $G$ is obtained from the disjoint union of a propeller or windmill and $k\ge 1$ copies of $P_3$ by identifying the nose of the propeller or windmill with exactly one leaf from each of the $P_3$s.

\end{itemize}

\newpage

\subsection{Connected, $K_4$-free, well-edge-dominated graphs with $7$ vertices}

\newpage
\subsection{Connected, $K_4$-free, well-edge-dominated graphs with $8$ vertices}

\end{document}